\documentclass[11pt]{article}
\usepackage{amssymb}
\usepackage{pstricks}
\usepackage{pst-node}
\begin{document}



\begin{center}
\textbf{THE EXISTENCE PROBLEM FOR STEINER NETWORKS IN STRICTLY CONVEX
DOMAINS}\\
by A. Freire\\
Knoxville, July 2006 \end{center} \vspace{.5cm}
\begin{center}ABSTRACT\end{center}
We consider the existence problem for `Steiner networks' (trivalent
graphs with $2\pi/3$ angles at each junction) in strictly convex
domains, with `Neumann' boundary conditions. For each of the three
possible combinatorial possibilities, sufficient conditions on the
domain are derived for existence; in addition, in each case explicit
examples of nonexistence are given. \vspace{.8cm}

\textbf{0. Introduction.}

About 20 years ago, motivated by dynamical models in materials
science describing phase separation and the motion of interfaces
separating phases, \cite{Bronsard-Reitich} introduced the problem
of motion of networks of curves in a planar domain with normal
velocity proportional to the curvature (at regular points), and
fixed angle conditions at the junctions. They derived (formally)
from the underlying model (an Allen-Cahn parabolic system with
two-dimensional order parameter and a three-well potential) the
geometric evolution, as well as the boundary conditions: the
angles formed by the curves at a `triple junction' (where three
arcs meet) is constant throughout the evolution, and the arcs
meeting the boundary of the domain do so orthogonally at all
times. In the same paper, short-time existence was proved for the
simplest network, three arcs meeting at an interior (moving)
point, with all angles equal to $2\pi/3$ radians. More recently,
\cite{Mantegazza et al.} undertook a thorough study of this system
of geometric evolution equations, which is distinguished by
non-standard boundary conditions at the junctions. In particular
(again for this simplest network, but also for fixed-endpoint
boundary conditions) they made explicit the natural continuation
criterion for the short-time solution, and developed the rescaling
analysis of singularity formation. This enabled them to rule out
finite-time curvature blowup at certain rates, which goes a long
way towards the natural global existence theorem in this setting.
 \vspace{.2cm}

With the goal of understanding the global non-linear evolution of
this system, it is natural to consider the existence,
classification and linearized stability properties of steady-state
solutions. (The flow is the gradient of total length of the
network, at least as long as the solution is classical, that is,
does not go through a topological change.) These have a very
simple variational description: consider embedded networks of
curves in a bounded domain $D\subset \mathbb R^2$, $C^1$ up to
each node and to $\partial D$, and with all nodes trivalent. If we
consider variations which preserve the combinatorics of the
network, but allow the boundary vertices to move freely on
$\partial D$, then the critical points of total length are exactly
those networks where (i) all edges are straight line segments;
(ii) the angles between edges at each node are $2\pi/3$; (iii) the
edges that meet $\partial D$ do so orthogonally. For brevity (and
by analogy with the fixed-boundary node case) we call these
`Steiner networks'. In a recent paper, \cite{Ikota-Yanagida}
settle the linearized stability question for Steiner networks that
are trees, without the assumption of convexity for the
domain.\vspace{.2cm}

It might seem at first that the existence of such networks of
arbitrary complexity, at least in the case of convex domains,
presents no greater subtlety. But already the first attempt at a
simple geometric construction to produce hexagonal cells in a
given domain led the author (a little over a year ago) to discover
a phenomenon akin to a simple kind of `holonomy' (\emph{section
5}). Once this case appeared settled, it was natural to try to
identify all the possible networks, and consider their existence
in a given convex domain. This is in principle a simple geometric
problem, but it leads to results that may seem surprising (at
least, they were surprising to the author.)\vspace{.2cm}

For example, consider complexity. One sometime sees optimistic
drawings of `honeycombs' of adjacent hexagonal cells tiling planar
domains. But, in fact, it is not too hard to show (\emph{Section
2}) that in a \emph{strictly} convex domain, there are only three
possible Steiner networks, excluding critical chords (Fig. 8): a
single triple junction, or `triode' (in the terminology of
\cite{Mantegazza et al.}); a single `double triode'; or a single
hexagonal cell, `anchored' to the boundary by six edges. In
addition, in such domains a Steiner network is necessarily
connected.\vspace{.2cm}

For the simplest case, the triode, one might hope (on general
`variational' grounds) that they always exist in smooth convex
domains. Indeed there have been such claims in the literature
(e.g. \cite{Tabachnikov}) based on the fact that the vertices of a
triode (when it exists) are contact points of a circumscribed
equilateral triangle which has critical perimeter among such
circumscribed triangles. But this ignores the fact that, for such
a `critical' configuration of boundary points to define a triode,
the inner normals at these points must meet \emph{inside $D$}. And
it is easy to construct examples of strictly convex domains in
which \emph{some} critical configurations violate
this.\vspace{.2cm}

Finding convex domains for which \emph{all} critical
configurations have the property that the normals meet outside $D$
proved harder; after spending some time trying to prove this
cannot happen, the author became aware of the class of convex
curves of `constant height', beautifully described in
\cite{Yaglom-Boltianskii}. For these curves the perimeter of a
circumscribed equilateral triangle is constant along the curve, so
the three inner normals at the contact points always meet at a
single point; and in one case, the \emph{constant-height biangle}
(which has two corners), the locus of these normal intersections
is entirely \emph{outside} the domain. This example can be
smoothed, and although the constant-height smoothings do support
triodes, the notion of `Minkowski sum' of convex domains suggests
a way to perturb them so as to destroy all triodes. Eventually it
turned out to be more efficient to implement all the constructions
completely analytically, using the support function in the
tangent-angle parametrization; this setup is described in
\emph{Section 1.} For the case of triodes, perhaps the most
striking non-existence result is Corollary 3.5: let $\cal C$ be a
smooth, strictly convex curve of constant height, with the
property that the three normals at contact points meet on the
outside, for some circumscribed equilateral triangle. \emph{Then
arbitrarily close to $\cal C$, one finds (i) strictly convex
curves which stably support triodes; (ii) strictly convex curves
which support no triodes, also stably.} (Here `arbitrarily close'
and `stably' refer to the strong $C^r$ topology in the
tangent-angle parametrization, for some $r\geq 2$.) On the
positive side, as expected, under a `pinching condition' for the
curvature, the normals always intersect internally, so at least
two triodes exist (Proposition 3.2).\vspace{.2cm}

The situation for double triodes (\emph{Section 4}) is somewhat
similar; in addition to a `criticality' condition and the
requirement that normals meet internally, a third condition plays
a role (note that the circle does not support double triodes.) For
non-existence results, curves of `constant width' exhibit a
similar `instability property' as constant-height curves for
triodes (Proposition 4.4). \emph{If ${\cal C}_0$ is a strictly
convex curve of constant width, arbitrarily close to ${\cal C}_0$
one finds: (i) strictly convex curves $\cal C$ for which all
sufficiently far outer parallel curve support double triodes;
(ii)strictly convex curves $\cal C$ for which neither $\cal C$
itself nor any outer parallel curve supports double triodes. } It
is harder to state a satisfactory sufficient criterion for
existence (the case of the circle rules out `curvature pinching').
Section 4 contains three existence results. Proposition 4.3 states
that \emph{if $\cal C$ is a strictly convex curve with only two
critical chords, making an angle greater than $\pi/3$, some outer
parallel curve stably supports a double triode.} In fact, if $\cal
C$ satisfies, in addition, the `curvature pinching' condition of
Prop. 3.2, $\cal C$ itself already supports double triodes.
\vspace{.2cm}

\emph{Section 5} deals with existence for hexagonal cells. Here
the criticality condition is vanishing of the `holonomy' , and one
has to consider a configuration of 24 points, at least some of
which must be inside the domain to guarantee existence. It turns
out that existence in a sufficiently far outer parallel curve of a
given strictly convex curve $\cal C$ always holds, and existence
of hexagonal cells for $\cal C$ itself is guaranteed under a
`curvature pinching' condition (Proposition 5.8). One gets
examples exhibiting neither hexagonal cells nor triodes by a
slight modification of the construction in section 3 (Corollary
5.10). \vspace{.2cm}

A natural question we are unable to resolve at this point is
whether there are strictly convex domains supporting no Steiner
networks at all (except for critical chords.) The problem is that
our construction of domains without triodes relies on the
existence of convex curves of constant height with exterior
intersections, and such examples are rare- we have been unable to
further `engineer' their width function so as to rule out double
triodes.\vspace{.2cm}

The results in the paper suggest that, also for the dynamical
problem, results obtained for evolution in the standard disk may
not carry over unchanged to more general convex domains; in
addition, we expect some of the techniques introduced here may be
useful to obtain estimates or constraints for the evolution
problem, at least in case the arcs remain convex. It is easy to
think of slightly more general settings in which some of the
phenomena should persist: for example, what if one removes the
requirement of \emph{strict} convexity? What about existence of
doubly-periodic Steiner networks (i.e., networks on a flat torus)?
Finally, one could look at the problem of geodesic Steiner
networks on surfaces of positive curvature. It is well-known that
those on the standard sphere have been classified, and the result
plays a role in the study of minimal surface sheets meeting at a
4-junction in the unit ball. If one looks at smooth ovaloids in
$\mathbb R^3$ (boundaries of smooth, strictly convex sets), it is
natural to wonder whether ovaloids exist which admit no geodesic
Steiner networks at all (closed geodesics don't
count.)\vspace{.2cm}

The study of Steiner networks with Dirichlet boundary conditions
(that is, critical-length graphs spanning a given set of points) has
a long history, both in the euclidean plane and on Riemannian
surfaces (see e.g. \cite{IvanovTuzhilin} for a survey.) For Steiner
networks of surfaces, a parametric version of the singular Plateau
problem was addressed recently via energy functionals in
\cite{MeseYamada}. Existence for the free-boundary case addressed in
the present paper (networks with endpoints moving freely on a given
boundary) does not seem to have been considered previously, with the
exceptions noted earlier.

This work was carried out entirely at the University of Tennessee,
Knoxville. It is a pleasure to thank my colleagues Nicholas Alikakos
(University of Athens), for his stimulating interest in this work
(which arose from a long-term joint project) during its development;
and Santiago Betel\'{u} (University of North Texas), for asking the
questions that got me thinking about existence for hexagonal
cells.\vspace{.4cm}


\newpage

\textbf{1. Strictly convex curves}. \vspace{.2cm}

\emph{1.1 Generalities.} We consider strictly convex oriented
Jordan curves $\cal C$, boundary of convex domains $D$ in $\mathbb
R^2.$ `Strictly convex' means that, for each $\theta$ in the unit
circle $S=\mathbb R\mbox{ mod }2\pi \mathbb Z$,  the oriented
support line of $D$ with unit normal $N(\theta)=(-\sin \theta,
\cos \theta)$ meets $\cal C$ at exactly one point. This defines a
continuous surjective map $B: S\rightarrow \cal C$, as well as a
continuous function $p(\theta)=-\langle
B(\theta),N(\theta)\rangle$. We are primarily interested in the
case where $B$ is injective and smooth (or at least piecewise
$C^2$), and thus defines a diffeomorphism $S\rightarrow \cal C$,
the `tangent angle parametrization'; for each $\theta \in S$, the
unit tangent vector to $\cal C$ at $B(\theta)$ is $T(\theta)=(\cos
\theta, \sin \theta)$, while $N(\theta)$ is the inner unit normal.
When $B$ is $C^1$, so is $p$ (the `support function'), and $B$ can
be recovered from $p$ via:
$$B(\theta)=-p(\theta)N(\theta)+p'(\theta)T(\theta).\quad(1.1)$$
If the origin $0\in \mathbb R^2$ is strictly inside $D$ (as will
be assumed throughout the paper), from $\langle
B(\theta),N(\theta)\rangle<0$ for all $\theta$ follows
$p(\theta)>0$. Since we only use the tangent angle
parametrization, relation (1.1) justifies describing classes of
curves, and indeed carrying out all of the analysis, in terms of
properties of $p$. If $p\in C^2(S)$, differentiating (1.1) we
obtain:
$$B'(\theta)=r(\theta)T(\theta),\quad \mbox{ where }
r(\theta)=p''(\theta)+p(\theta).$$ Thus, where $r(\theta)>0$, $B$
is an immersion with unit tangent vector $T(\theta)$;
geometrically, $r(\theta)$ is the radius of curvature at
$B(\theta)$. Note that the tangent-angle parametrization is less
regular than arc length: if the arc length parametrization
$\Gamma(s)$ is $C^2$ with positive curvature
$k=\frac{d\theta}{ds}=\langle
\Gamma_{ss},\Gamma_s^{\perp}\rangle$, we have $r=1/k$ positive and
continuous, so $B$ is only $C^1$. The support function $p$ can be
recovered from $r$ by integration (assuming $r$, say, piecewise
$C^0$:)
$$p(\theta)=p(0)\cos(\theta)+p'(0)\sin(\theta)+
\int_0^{\theta}\sin(\theta-\tau)r(\tau)d\tau.\quad (1.2)$$ With
$r>0$ of class $C^k$ and $2\pi$-periodic, (1.1) and (1.2) will
always define a parametrization of class $C^{k+2}$ of a strictly
convex Jordan curve, provided only $p$ given by (1.2) is
$2\pi$-periodic, for which it suffices to impose:
$$\int_0^{2\pi}r(\tau)\cos(\tau)d\tau=\int_0^{2\pi}r(\tau)\sin
\tau d\tau=0.$$ (This is automatic if $r=p''+p$ with $p$
$2\pi$-periodic and piecewise $C^2$.) Our main interest is in
smooth (or $C^r$ for some $r\geq 2$) strictly convex curves,
identified with the set:
$${\cal P}_{smooth}=\{p\in C^2(\mathbb R), 2\pi\mbox{-periodic};
p>0, r=p''+p>0\}.$$ We will also need to consider piecewise $C^2$
strictly convex curves without corners:
$${\cal P}_{pw}=\{p\in C^2_{pw}(\mathbb R), 2\pi\mbox{-periodic};
p>0, r=p''+p>0\},$$ as well as piecewise $C^2$ strictly convex
curves with corners:
$${\cal P}_{corner}=\{p\in C^2_{pw}(\mathbb R),
2\pi\mbox{-periodic}; p>0, r=p''+p\geq 0,
Z(r)=\bigcup_{i=1}^N[a_i,b_i]\},$$ where the last condition means
the zero set $Z(r)$ of $r$ in $[0,2\pi]$ is assumed to be a finite
union of closed non-degenerate intervals.\vspace{.2cm}

If $p\in {\cal P}_{corner}$, let $I=[a,b]$ be a $\theta$ interval
on which $r=0$. Then $B'(\theta)=0$ for $\theta \in I$, and
$B\equiv B(a)$ in $I$. Thus $B$ is not, strictly speaking, a
parametrization on the interval $I$, and the unit normals
$\{N(\theta); \theta \in I\}$ correspond to the set of support
lines at $B(a)$, the `corner' corresponding to $I$. Adding a
positive constant $c$ to an element $p$ of ${\cal P}_{corner}$
yields an element of ${\cal P}_{pw}$, that is, erases the corners.
Geometrically, this corresponds to considering the `exterior
parallel curves' $B_c(\theta)$ to the convex curve defined by $p$:
from (1.1), we see that:
$$B_c(\theta)=-(p+c)N(\theta)+p'(\theta)=B(\theta)-cN(\theta).$$
More generally, it is clear that each of the three classes ${\cal
P}_{smooth}\subset{\cal P}_{pw}\subset{\cal P}_{corner}$ is closed
under linear combinations with positive coefficients. This
corresponds to the well-known fact that convexity is preserved by
homothety (with a fixed center $0$) and Minkowski sum; the latter
operation is conveniently described in terms of the tangent-angle
parametrization (or `pseudo-parametrization', in the case of ${\cal
P}_{corner})$: if $B_1(\theta)$, $B_2(\theta)$ parametrize ${\cal
C}_1$ (resp. ${\cal C}_2$), the Minkowski sum $D_1+D_2$ of the
domains they bound has support function $p_1(\theta)+p_2(\theta)$,
and its boundary is (pseudo)parametrized by
$B_1(\theta)+B_2(\theta)$. Adding $c>0$ to $p$ corresponds to taking
Minkowski sum with a disk of radius $c$. \vspace{.2cm}

We can use the support function in the tangent-angle
parametrization to introduce a strong topology in the space of
strictly convex curves: the uniform $C^2$ topology induced in the
convex cone ${\cal P}_{smooth}$ from $C^2(S)$. Accordingly, we say
something happens `stably' for a given $p_0\in {\cal P}_{smooth}$
(or even $p_0\in {\cal P}_{corner}$) if it also happens for all
$p\in {\cal P}_{smooth}$ in some $C^2$-neighborhood of $p_0$.

\vspace{.3cm}

\emph{1.2 Special classes of convex curves.} \vspace{.2cm}

\emph{1.2.1 Width and constant width.} Let $\cal C$ be a strictly
convex curve in $\mathbb R^2$. The \emph{width} $b(\theta)$ in
direction $\theta \in S$ is the unoriented distance between the
supporting lines with unit normals $N(\theta), N(\theta+\pi)$.
Clearly, since $0\in D$:
$$b(\theta)=p(\theta)+p(\theta+\pi)=\langle
B(\theta+\pi)-B(\theta),N(\theta)\rangle>0.$$ Width behaves
linearly under homothety/Minkowski sum (at corresponding
directions). Denote by $w\in C^1_{pw}$ the `derived width':
$$w(\theta)=-b'(\theta)=-p'(\theta)-p'(\theta+\pi)=\langle
B(\theta+\pi)-B(\theta),T(\theta)\rangle.$$

We say $\cal C$ is \emph{symmetric} if $B(\theta+\pi)=-B(\theta)$
for all $\theta$; this is equivalent to $\pi$-periodicity for the
support function $p(\theta)$, so in this case
$b(\theta)=2p(\theta)$. It is easy to see that one can always
construct symmetric convex curves with prescribed `derived width':
\vspace{.2cm}

\textbf{Proposition 1.1.} \emph{Let $w$ be $C^1$, $\pi$-periodic
and satisfy $\int_0^{\pi}w(\tau)d\tau=0$. There exists a
one-parameter family of `parallel' strictly convex symmetric
curves with `derived width function' $w$.}

\emph{Proof.} Define $b(\theta)=C-\int_0^{\theta}w$; $b$ is $C^2$,
$\pi$-periodic, positive if $\max_{\theta}\int_0^{\theta}w<C$ and
satisfies $b''+b>0$ if $\max_{\theta}\{\int_0^{\theta}w(\tau)d\tau
+w'(\theta)\}<C$. Both conditions will hold for $C$ in some
interval $(C_0,\infty)$, and then $p(\theta)=(1/2)b(\theta)\in
{\cal P}_{smooth}$ and defines a strictly convex symmetric curve
with `derived width' $w$. \vspace{.2cm}

\emph{Example 1.1.} The \emph{Reuleaux triangle} of constant width
$b$ is obtained from an equilateral triangle $B_1B_2B_3$ of side
length $b$ by drawing three circular arcs with center $B_i$,
radius $b$ and aperture $\pi/3$ radians (Fig. 1). Placing the
origin of $\mathbb R^2$ at the barycenter of the triangle, and
fixing $B(0)=(p'(0),-p(0))=(0,-b(1-\frac{\sqrt 3}3))$ gives the
support function $p\in {\cal P}_{corner}$:
$$p(\theta)=\left \{ \begin{array}{ll}
&b-b\frac{\sqrt{3}}3\cos \theta, \theta \in [-\pi/6, \pi/6]\\
&\frac b2\sin \theta +b\frac{\sqrt{3}}6\cos \theta, \theta \in
[\pi/6,\pi/2],\end{array}\right .\quad (1.3)$$ extended to
$\mathbb R$ with period $2\pi/3$.
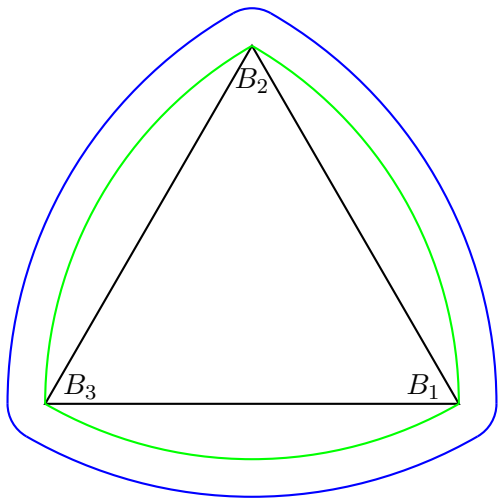
\begin{figure} \psset{unit=0.5cm}
\begin{center}
\begin{pspicture}(-7,-13)(7,1)
 \SpecialCoor
\pnode(0,0){A}\pnode(11;240){B}\pnode(11;300){C}
\pspolygon(A)(B)(C)
\psarc[linecolor=green](A){11}{240}{300}
\psarc[linecolor=green](B){11}{0}{60}
\psarc[linecolor=green](C){11}{120}{180}
\psarc[linecolor=blue](A){12}{240}{300}
\psarc[linecolor=blue](B){12}{0}{60}
\psarc[linecolor=blue](C){12}{120}{180}
\psarc[linecolor=blue](A){1}{60}{120}
\psarc[linecolor=blue](B){1}{180}{240}
\psarc[linecolor=blue](C){1}{300}{0}
\uput{0.6}[270](A){$B_2$}
\uput{0.5}[30](B){$B_3$}\uput{0.5}[150](C){$B_1$}
\end{pspicture}
 
\caption{ The Reuleaux triangle and a smoothed parallel
curve}\end{center}
\end{figure}

The arc of the curve corresponding to $\theta \in [-\pi/6,\pi/6]$
is regular, while $B(\theta)\equiv B_1= b(1/2,-\sqrt{3}/6)$ (a
vertex) for $\theta \in [\pi/6, \pi/2]$. It is geometrically clear
that $p$ is $2\pi/3$-periodic, and therefore
$B(\theta+2\pi/3)=R_{2\pi/3}B(\theta)$ for all $\theta$ ($\cal C$
is `3-symmetric', but not symmetric; the only symmetric curve of
constant width is the circle.) One checks easily that $p\in C^1$,
and $B(\theta)$ maps $[\pi/2, 5\pi/6]$ and $[\pi/6, 3\pi/2]$ to
regular arcs and $[5\pi/6, 7\pi/6], [3\pi/2, 11\pi/6]$ to the
vertices $B_2=(0, b\sqrt{3}/3)$ and $B_3=b(-1/2,-\sqrt{3}/6)$.

The radius of curvature $r(\theta)$ is constant (equal to $b$) on
the regular intervals, and vanishes on the singular intervals.
Adding an arbitrary positive constant $c$ to $p$ yields a support
function in ${\cal P}_{pw}$; the corresponding $B_c(\theta)$
parametrizes the $c$-parallel curve, which is of constant width
$b+c$. Although this curve has no corners, the parametrization
$B_c(\theta)$ is only piecewise $C^1$. \vspace{.2cm}

It will be of interest to consider explicit examples of curves of
constant width with smooth support function, and it is natural to
use truncated Fourier series. In general, write:
$$p(\theta)=C+\sum_{n\geq
1}a_n\cos{n\theta}+b_n\sin{n\theta},$$ where $C,a_n,b_n$ are
Fourier coefficients. We have the following simple observations:

(i) $\cal C$ is centrally symmetric ($p(\theta)=p(\theta+\pi),
B(\theta+\pi)=-B(\theta)$) iff $a_n$ and $b_n$ vanish for $n$
\emph{odd};

(ii) $\cal C$ is 3-symmetric ($p$ is $2\pi/3$-periodic) iff $a_n$
and $b_n$ vanish unless $n$ is a \emph{multiple of 3;}

(iii) $\cal C$ has constant width iff $a_n$ and $b_n$ vanish for
$n$ \emph{even}.

(iv) Suppose the $x$-axis is an axis of symmetry of $\cal C$. With
$B(\theta)=(b_1(\theta),b_2(\theta))$, this means
$b_2(\theta)+b_2(\pi-\theta)\equiv 0$, or equivalently
$p(\pi-\theta)=p(\theta)$ for all $\theta$. In terms of Fourier
coefficients, this means: $a_n=0$, $n$ odd; $b_n=0$, $n$ even. (In
this case the $x$- axis is a critical chord, and there is a second
critical chord perpendicular to it.)

(v) Now suppose the $y$-axis is an axis of symmetry of $\cal C$.
Then $b_1(\theta)+b_1(2\pi-\theta)\equiv 0$. Equivalently, $p$ is
even in $\theta$, so $b_n=0$ for all $n$. \vspace{.2cm}

\emph{Example 1.1, cont.} Thus the Fourier coefficients for the
Reuleaux triangle vanish unless $n$ is an \emph{odd multiple of 3}
($n=3,9,15,\ldots$), and all sine coefficients $b_n$ vanish. For
the first few nonzero terms we find:
$$ p(\theta)=\frac 12 -\frac 1{4\pi}\cos{3\theta}+\frac
1{120\pi}\cos{9\theta}-\frac 1{560\pi}\cos{15\theta}+\ldots$$ We may
obtain a smooth curve of constant width by taking a finite number of
terms; but note that, since $r(\theta)$ is only piecewise
continuous, the radius of curvature for the approximation may fail
to be positive, no matter how many terms we take (Gibbs phenomenon);
so we have to add a constant to $p$. For example, taking the first
three terms given above we find that the minimum curvature (attained
at $\theta=\pi/4$) is $\frac
12-\frac{4\sqrt{2}}{3\pi}=-0.1002\ldots$, and subtracting this
number from $p$ we obtain:
$$p_{\epsilon}(\theta)=\frac{4\sqrt{2}}{3\pi}-\frac
1{4\pi}\cos{3\theta}+\frac 1{120\pi}\cos{9\theta}+\epsilon.\quad
(1.4)$$ For any $\epsilon >0$, $p_{\epsilon}\in {\cal P}_{smooth}$
has constant width, and in addition is 3-symmetric and symmetric
with respect to the $y$-axis; the corresponding $B_{\epsilon}$
parametrizes a `smoothed Reuleaux triangle' (shown in
Fig.1).\vspace{.3cm}

\emph{1.2.2. Height and constant height.}

Given an oriented convex curve $\cal C$ and a direction $\theta$,
there is a unique smallest oriented equilateral triangle $\cal
T(\theta)$ enclosing $\cal C$ with sides parallel to $T(\theta)$,
$T(\theta+2\pi/3), T(\theta+4\pi/3)$. If $\cal C$ is strictly
convex, possibly with corners, $\cal T(\theta)$ has exactly three
points in common with $\cal C$: $B_1=B(\theta), B_2=B(\theta
+2\pi/3), B_3=B(\theta+4\pi/3)$; each side of the triangle is
contained in a supporting line.

Referring to Fig.9 (section 3), we have for the vertices of $\cal
T(\theta)$:
$$Q_1=B_1+\frac 2{\sqrt{3}}s_1T_1=B_2+\frac 2{\sqrt{3}}t_2T_2.$$
Taking inner products with $T_1$ and with $T_2$, we find:
$$s_1=-\langle B_2-B_1,N_2 \rangle,\quad t_2=\langle B_1-B_2,N_1\rangle,$$
and cyclically mod 3:
$$Q_3=B_1+\frac 2{\sqrt{3}}t_1T_1, \mbox{ with }t_1=\langle B_1-B_2,
N_1\rangle. $$ Thus the side length of $\cal T(\theta)$ is given
by: $$\begin{array}{ccc} l_{tri}(\theta)=\langle
Q_1-Q_3,T_1\rangle =\frac 2{\sqrt{3}}(s_1-t_1)&=&-\frac
2{\sqrt{3}}(\langle
B_2-B_1,N_2\rangle+\langle B_3-B_1, N_3\rangle)\\
&=& -\frac 2{\sqrt{3}}(\langle B_2, N_2\rangle+\langle
B_3,N_3\rangle+\langle B_1,N_1\rangle),\end{array}$$ using
$N_1+N_2+N_3=0$. We use the notation $N_i=N(\theta+2(i-1)\pi/3)$
throughout the paper, with $i=1,2,3$.

Thus we see that the height of the triangle $\cal T(\theta)$ is:
$$h(\theta)=p(\theta)+p(\theta+2\pi/3)+p(\theta+4\pi/3).$$
We call $h(\theta)$ the \emph{height function} of $\cal C$ (or of
$p$) and define the `derived height' by:
$$\begin{array}{ccc}
\omega(\theta)=h'(\theta)&=&\langle B_1,T_1\rangle +\langle B_2,
T_2\rangle + \langle B_3, T_3\rangle\\
&=&\langle B_2-B_1,T_2\rangle-\langle
B_1-B_3,T_3\rangle,\end{array}\quad (1.5)$$ since $T_1+T_2+T_3=0$.
\vspace{.2cm}

In complete analogy with Proposition 1.1, one can always find
strictly convex curves with given `derived height'. The proof is
completely analogous.

\textbf{Proposition 1.2.} \emph{Let $\omega$ be $C^1$ and
$2\pi/3$-periodic, and satisfy
$\int_0^{2\pi/3}\omega(\tau)d\tau=0$. There exists a one-parameter
family of parallel, 3-symmetric strictly convex curves with
derived height function $\omega$.} \vspace{.2cm}

One checks easily that to the list of characterizations of
properties of $\cal C$ by the Fourier coefficients of $p$ one may
add:

(vi) $\cal C$ is a curve of constant height ($\omega\equiv 0$) iff
$a_n=b_n=0$ whenever $n$ is a multiple of 3. \vspace{.2cm}

The only 3-symmetric curve of constant height is therefore the
circle. \vspace{.2cm}

\emph{Example 1.2.} There are also strictly convex curves of
constant height with corners, consisting of any number of circular
arcs (except for multiples of 3) with the same radius (see
\cite{Yaglom-Boltianskii}.) We describe in detail the simplest of
them, the \emph{constant height biangle} of height $h$.

Starting from the equilateral triangle $B_1B_2B_3$ of height $h$
(with $B_1-B_3=hT(0)$, see Fig.2), we draw a $\pi/3$ arc with
center $B_2$ and radius $h$ (which intersects the triangle at
$M,N$), then reflect this arc on the segment $MN$. Taking as the
origin of $\mathbb R^2$ the midpoint of $MN$, we compute the
support function $p$:
$$p(\theta)=\left \{ \begin{array}{ll}&h(1-\frac {\sqrt{3}}2, \cos
\theta), \theta \in [-\pi/6, \pi/6],\\
&\frac h2 \sin \theta, \theta \in [\pi/6,
5\pi/6],\end{array}\right .\quad(1.6)$$ extended to the real line
as a $\pi$-periodic function. One checks easily that $p\in {\cal
P}_{corner}$ (in particular, $p\in C^1$), and the corresponding
$B(\theta)$ maps $[-\pi/6,\pi/6]$ and $[5\pi/6, 7\pi/6]$ to the
regular arcs $MB(0)N$, $NB(\pi)M$ (resp.) and $[\pi/6,
5\pi/6],[7\pi/6,11\pi/6]$ to the vertices
$N=(\frac{h\sqrt{3}}2,0), M=(-\frac{h\sqrt{3}}2,0)$(resp.) Adding
positive constants to $p$ we obtain the support functions of the
outer parallel curves of the biangle, which are constant-height
curves without corners (but only piecewise $C^2$).

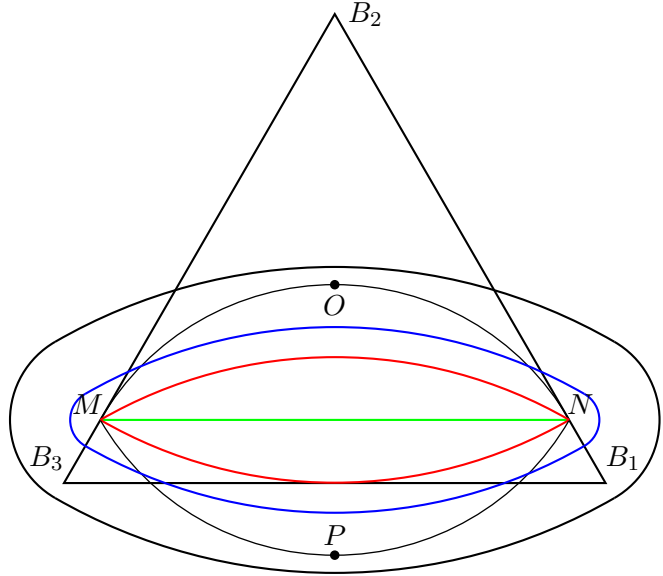
\begin{figure}\psset{unit=0.8cm}
\begin{center}

\begin{pspicture}(-6,-5)(6,5)
\psset{unit=0.8cm}  \SpecialCoor
\pnode(3.3;270){H} \pnode([angle=180,nodesep=4.5]H){Q30}
\pnode([angle=0,nodesep=4.5]H){Q10}\pnode(4.5;90){Q20}
\pspolygon(Q10)(Q20)(Q30) \pnode([angle=240,nodesep=7.794]Q20){M}
\pnode([angle=300,nodesep=7.794]Q20){N}
\psline[linecolor=green](M)(N)
\psarc[linecolor=red](Q20){7.794}{240}{300} \pnode(4.5;270){P0}
\psarc[linewidth=0.5pt](P0){4.5}{30}{150}
\psarc[linewidth=0.5pt](0,0){4.5}{210}{330} \pnode(9;270){Q2bar}
\psarc[linecolor=red](Q2bar){7.794}{60}{120}
\psarc[linecolor=blue](Q20){8.294}{240}{300}
\psarc[linecolor=blue](Q2bar){8.294}{60}{120}
\psarc[linecolor=blue](M){0.5}{120}{240}
\psarc[linecolor=blue](N){0.5}{300}{60}
\psarc(Q20){9.294}{240}{300} \psarc(Q2bar){9.294}{60}{120}
\psarc(M){1.5}{120}{240} \psarc(N){1.5}{300}{60}
\uput{0.1}[120](M){$M$}\uput{0.1}[60](N){$N$}
\uput{4}[90](P0){$O$}\uput{4}[270](0,0){$P$}
\uput[60](Q10){$B_1$}\uput[0](Q20){$B_2$}\uput[120](Q30){$B_3$}
\psdots(4.5;270)([angle=90,nodesep=4.5]P0)
\end{pspicture}
 
\end{center}
\caption{The constant-height biangle and two smoothed parallel
curves; also shown is the locus $OMPN$ of normal intersections for
all three curves.}
\end{figure}

As before, to obtain smooth examples we consider Fourier series.
Since the biangle is symmetric with respect to the $x$ and $y$
axes, and of constant height, from (iv), (v) and (vi) above we see
that its support function has a Fourier cosine series, with $a_n$
non-vanishing only for $n$ even, and not a multiple of 3. The
first few terms are (we set $h=1$):
$$p(\theta)=\frac
13-\frac{\sqrt{3}}{3\pi}\cos(2\theta)-\frac{\sqrt{3}}{30\pi}\cos(4\theta)
+\frac{\sqrt{3}}{252\pi}\cos(8\theta)+\ldots$$ The radius of
curvature corresponding to the first three terms is:
$$r(\theta,4)=\frac 13+\frac{\sqrt{3}}\pi (\cos(2\theta)+\frac 12 \cos
(4\theta)),$$ attaining the minimum value $\frac
13-\frac{3\sqrt{3}}{4\pi}=-0.0802\ldots$ at $\theta=\pi/3$. We
subtract this value from $p$ to obtain the support function of a
`smoothed biangle':
$$p_{\epsilon}(\theta)=\frac{\sqrt{3}}\pi(\frac 34-\frac
13\cos(2\theta)-\frac 1{30}\cos(4\theta))+\epsilon.\quad (1.7)$$
For any $\epsilon >0$, this defines a support function in ${\cal
P}_{smooth}$, and the corresponding curve is strictly convex, of
constant height, and has the symmetries of the biangle.(See
Fig.2).

The constant-height biangle and its smoothed versions will play an
important role in the construction of examples.


\newpage

\textbf{2. Steiner networks in strictly convex domains.}
\vspace{.2cm}

A \emph{network} in a bounded convex domain $D\subset \mathbb R^2$
is an embedded (undirected) graph $\cal N$ intersecting ${\cal
C}=\partial D$ at finitely many points (all univalent vertices of
$\cal N$), and with all interior vertices \emph{trivalent.} We
assume each edge admits a regular $C^1$ parametrization up to its
end-vertices and/or the boundary; in particular, the total length
$L$ of $\cal N$ is well-defined. The critical points of $L$ (with
respect to variations that don't change $\cal N$ combinatorially,
but allow the boundary vertices to move on the boundary) are
`\emph{Steiner networks}' with free (Neumann) boundary conditions.
By definition, this means the edges are line segments, meeting at
each interior vertex with three angles equal to $2\pi/3$ radians,
and meeting $\partial D={\cal C}$ orthogonally at the boundary
vertices. The main problem addressed in this paper is the
existence and classification of such networks in strictly convex
planar domains. (A special case is a `critical chord' in $D$, a
one-edge network with no interior vertices.)\vspace{.2cm}

In this section we show that the combinatorial possibilities for
such networks in a strictly convex domain are in fact rather
limited-other than critical chords, only three different types may
occur (Fig.8): the triple junction or `triode' (in the terminology
of \cite{Mantegazza et al.}); the `double triode'; and a single
`hexagonal cell', anchored to the boundary by six edges.
Furthermore, a Steiner network is necessarily connected, although
this is not assumed \emph{a priori.}

In the following, $\cal N$ denotes a Steiner network in a strictly
convex bounded domain $D\subset \mathbb{R}^2$. \vspace{.2cm}

\emph{2.1 Chains.} A \emph{chain} $C$ is a connected subgraph of
$\cal N$, including at least two edges, without branching (each
internal vertex of $C$ is adjacent to two edges of $C$) and such
that all exterior angles between two consecutive edges (in
principle either $\pi/3$ or $-\pi/3$) have the same sign. We may
choose a consistent orientation of the edges of $C$ so that all
exterior angles are equal to $\pi/3$; the chain then has an
initial vertex $v_0$ and a final vertex $v_N$. A chain is
`inextendible' if it is not a subset of another chain. Any chain
can be continued to an inextendible one, in a unique way.
\vspace{.2cm}

(1) Claim: \emph{An inextendible chain either starts and ends at
two different boundary vertices of $\cal N$ or is closed
($v_N=v_0$), and then consists of six edges and six interior
vertices, making up a convex equiangular hexagon.}\vspace{.2cm}

\emph{Proof.} Assume, by contradiction, $v_N$ is an interior
vertex and $C$ is not closed. Since the chain is not
forward-extendible past $v_N$, $v_N$ must be adjacent (in $\cal
N$) to another vertex $v$ of $C$. $v_N$, $v$, the edge $f$
connecting them and the edge $e$ of $C$ arriving at $v$ together
define a configuration as shown in Fig. 3.\vspace{.2cm}

From the definition of `chain', there is an arc of $C$ (i.e., a
subchain $C_1$) which (together with the edge $f$ of $\cal N$
connecting $v_N$ and $v$) bounds an open region $U$ containing an
edge $e$ of $C$ arriving at $v$ (Fig. 3). \vspace{.2cm}

\begin{figure} \psset{unit=0.5cm}
\begin{center}

\SpecialCoor
\begin{pspicture}(-5,-6)(5,4)
\psline{->}(0,0)(2;210) \pnode(2;210){v1}
\psline{->}(v1)([angle=270,nodesep=3]v1)
\psline[linestyle=dashed](v1)([angle=150,nodesep=1.5]v1)
\pnode([angle=270,nodesep=3]v1){v2}
\psline{->}(v2)([angle=330,nodesep=2.7]v2)
\pnode([angle=330,nodesep=2.7]v2){v3}
\psline{->}(v3)([angle=30,nodesep=4.5]v3)
\pnode([angle=30,nodesep=4.5]v3){v4}
\psline[linestyle=dashed](v2)([angle=210,nodesep=3]v2)
\psline[linestyle=dashed](v4)([angle=330,nodesep=2]v4)
\psline{->}(v4)([angle=90,nodesep=8]v4)
\pnode([angle=90,nodesep=8]v4){v5}
\psline{->}(v5)([angle=150,nodesep=3.7]v5)
\psline[linestyle=dashed](v5)([angle=30,nodesep=2.2]v5)
\pnode([angle=150,nodesep=3.7]v5){v6}
\psline{->}(v6)([angle=210,nodesep=1.5]v6)
\pnode([angle=210,nodesep=1.5]v6){v7}
\psline[linestyle=dashed](v7)([angle=150,nodesep=1.5]v7)
\pnode(0,0){v} \psline{->}([angle=330,nodesep=1.5]v)(v)
\psline[linestyle=dashed](v)([angle=90,nodesep=5.9]v)
\psline[linestyle=dashed](v6)([angle=90,nodesep=1.5]v6)
\psline[linestyle=dashed](v3)([angle=270,nodesep=1.5]v3)
\psline[linestyle=dashed](v4)([angle=330,nodesep=1.5]v4)
\uput{0.5}[180](v){$v$} \uput{0.4}[330](v7){$v_N$}
\uput{1.2}[350](v){$e$}
\uput{3.4}[100](v){$f$}\uput{2}[45](v2){$U$}\uput{4}[80](v4){$C_1$}
\end{pspicture}
 
\caption{proof of claim 1, first part }\end{center}
\end{figure}
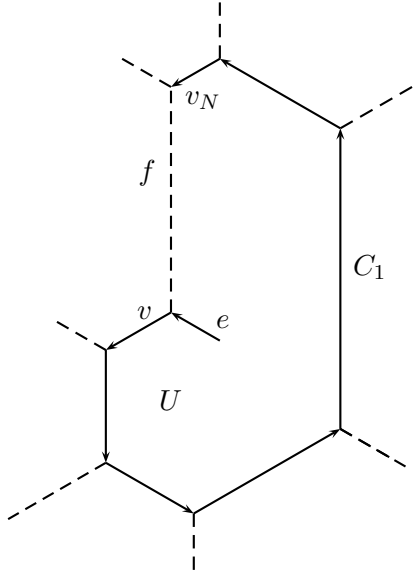

The backward continuation of the chain from $e$ is contained in
$U$, and does not meet $C_1$; in particular, it never meets the
boundary of $D$, so the initial vertex $v_0$ is also interior.
Since $C$ cannot be backward-continued further past $v_0$, $v_0$
is adjacent in $\cal N$ to a vertex $\bar{v}$ of $C$ preceding $v$
in $C$, and we must have the following configuration including
$v_0$, $\bar{v}$ and two edges $g,\bar{e}$ as shown in Fig.4.
\vspace{.2cm}

But then there must be a second sub-arc $C_2$ of $C$ (contained in
$U$) from $v_0$ to $\bar{v}$, which (together with the edge $g$ of
$\cal N$ connecting $v_0$ and $\bar{v}$) bounds an open set
$V\subset U$ containing the edge $\bar e$ of $C$ leaving $\bar v$
(Fig. 4). Forward continuation of $C_2$ past $\bar{e}$ can never
meet a vertex of $C_2$, hence is entirely contained in the open
set $V$, and will not meet the arc $C_1$ of the chain- a
\emph{contradiction.}

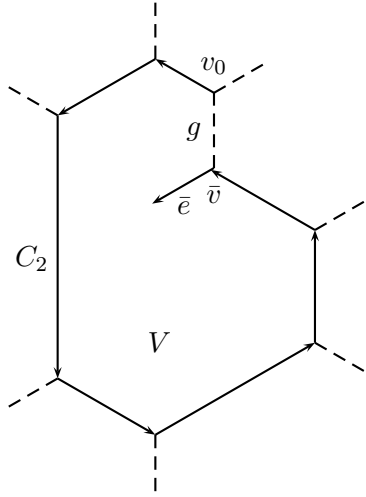
\begin{figure} \psset{unit=0.5cm}
\begin{center}

\SpecialCoor
\begin{pspicture}(-5,-8)(5,5)
\pnode(0,0){vbar} \pnode([angle=90,nodesep=2]vbar){v0}
\psline{->}(v0)([angle=150,nodesep=1.8]v0)
\pnode([angle=150,nodesep=1.8]v0){v1}
\psline{->}(v1)([angle=210,nodesep=3]v1)
\pnode([angle=210,nodesep=3]v1){v2}
\psline{->}(v2)([angle=270,nodesep=7]v2)
\pnode([angle=270,nodesep=7]v2){v3}
\psline{->}(v3)([angle=330,nodesep=3]v3)
\pnode([angle=330,nodesep=3]v3){v4}
\psline{->}(v4)([angle=30,nodesep=4.9]v4)
\pnode([angle=30,nodesep=4.9]v4){v5}
\psline{->}(v5)([angle=90,nodesep=3]v5)
\pnode([angle=90,nodesep=3]v5){v6}
\psline{->}(v6)([angle=150,nodesep=3.2]v6)
\psline[linestyle=dashed](vbar)(v0)
\psline{->}(vbar)([angle=210,nodesep=1.9]vbar)
\psline[linestyle=dashed](v0)([angle=30,nodesep=1.5]v0)
\psline[linestyle=dashed](v1)([angle=90,nodesep=1.5]v1)
\psline[linestyle=dashed](v2)([angle=150,nodesep=1.5]v2)
\psline[linestyle=dashed](v3)([angle=210,nodesep=1.5]v3)
\psline[linestyle=dashed](v4)([angle=270,nodesep=1.5]v4)
\psline[linestyle=dashed](v5)([angle=330,nodesep=1.5]v5)
\psline[linestyle=dashed](v6)([angle=30,nodesep=1.5]v6)
\uput{0.5}[270](vbar){$\bar{v}$} \uput{0.5}[90](v0){$v_0$}
\uput{0.8}[120](vbar){$g$}\uput{1}[230](vbar){$\bar{e}$}
\uput{2.5}[20](v3){$V$}\uput{3.5}[260](v2){$C_2$}
\end{pspicture}
 
\caption{proof of claim 1, second part }\end{center}
\end{figure}

The case when the contradiction hypothesis is that $v_0$ is
interior is completely analogous (or just reverse the
orientations.)\vspace{.2cm}

If the chain is closed, since each interior angle is $2\pi/3$, it
must be an equiangular hexagon. \vspace{.2cm}

(2) \emph{If $D$ is strictly convex, a maximal chain connecting
two boundary vertices cannot contain more than three edges.}

A fourth edge would give a total `turning angle' of at least $\pi$
for the chain, which is impossible for a chain connecting two
boundary points (by strict convexity). Note that (1) and (2) imply
that \emph{a chain including a boundary vertex has at most two
interior vertices} (since forward continuation from a third
interior vertex would yield an inextendible chain with at least
four edges, beginning and ending at boundary vertices.)
\vspace{.2cm}

\emph{2.2 Classification of connected components of $\cal N$.} If
we exclude the case of critical chords, each boundary vertex of
$\cal N$ is adjacent to a unique vertex, necessarily an interior
one. Let $b_0$ be a boundary vertex of a connected component
$\hat{\cal N}$, adjacent to the interior vertex $v_1$. There are
only three (unoriented) directions for edges in the network, so by
rotating $D$ we may assume $b_0$ is `vertically above' $v_1$, and
then phrases such as `upper left', `lower right', `vertically
below' have a well-defined meaning. Three cases are possible.

(i) If $v_1$ connects to two other boundary vertices $b_1, b_2$,
the connected component $\hat{\cal N}$ is a triode (Fig.5).

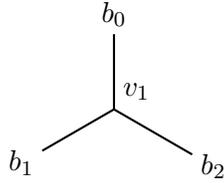
\begin{figure} \psset{unit=0.2cm}
\begin{center}
\SpecialCoor
\begin{pspicture}(-5,-5)(5,5)
 \psline(0,0)(5;90) \psline(0,0)(5.5;210)
\psline(0,0)(6;330)
\uput[ur](0,0){$v_1$} \uput{5.5}[90](0,0){$b_0$}
\uput{6}[210](0,0){$b_1$} \uput{6.5}[330](0,0){$b_2$}
\end{pspicture}
 
\caption{Classification, part(i)-triode}\end{center}
\end{figure}

(ii) (Fig.6) If exactly one of the vertices (other than $b_0$)
adjacent to $v_1$ in $\cal N$ is interior (say, the lower-right
vertex $v_2$, while the lower-left vertex $b_1$ is a boundary
vertex), both the remaining vertices adjacent to $v_2$ must be
boundary vertices; otherwise, there would be a chain beginning at
a boundary vertex ($b_0$ or $b_1$) and containing at least three
interior vertices. Thus, in this case the connected component
$\hat{\cal N}$ is a `double triode'.

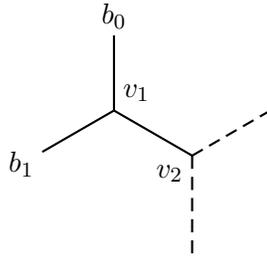
\begin{figure} \psset{unit=0.2cm}
\begin{center}
\begin{pspicture}(-8,-7)(8,5)\SpecialCoor
\psline(0,0)(5;90) \psline(0,0)(5.5;210) \psline(0,0)(6;330)
\pnode(6;330){A}
\psline[linestyle=dashed](6;330)([angle=30,nodesep=6]A)
\psline[linestyle=dashed](A)([angle=270,nodesep=6.5]A)
\uput[ur](0,0){$v_1$} \uput{5.5}[90](0,0){$b_0$}
\uput[dl](A){$v_2$} \uput{6}[210](0,0){$b_1$}
\end{pspicture}
 
\caption{Classification, part(ii)-double triode}\end{center}
\end{figure}

(iii) Assume now both vertices (other than $b_0$) adjacent to
$v_1$ in $\cal N$ are interior ($v_L$ to the left, $v_R$ to the
right; Fig.7). The remaining `upper' adjacent vertices of $v_L$
and $v_R$ are both boundary vertices ($b_L$ and $b_R$,
respectively), to avoid `long' boundary chains starting at $b_0$.
The `lower' adjacent vertices $v_L'$, $v_R'$ must be
\emph{interior} vertices; otherwise one would have `long' chains
connecting a boundary vertex ($v_L'$ or $v_R'$) to an interior
vertex ($v_R$ or $v_L$, resp.) Of the two lower adjacent vertices
to $v_L'$ (resp. $v_R'$) the one on the left (resp. right) must be
a \emph{boundary} vertex ($b_L'$, resp. $b_R'$), to avoid `long'
chains connecting the boundary vertex $b_L$ (resp. $b_R$) to an
interior vertex. The remaining adjacent vertices of $v_L'$, $v_R'$
must be interior (were  either of them a boundary vertex, it would
be part of a `long' chain including $v_1$); call them
$v_2^L,v_2^R$. We claim that, in fact, we must have
$v_2^L=v_2^R:=v_2$; equivalently, the chain
$C=v_2^Rv_R'v_Rv_1v_Lv_L'v_2^L$ is \emph{closed}.

Consider the inextendible chain $\hat{C}$ containing $C$. If $C$
is not closed, $\hat C$ must begin and end on the boundary of $D$,
and, as seen in (2) above, in this case $\hat{C}$ has at most
three edges. Since the chain $C$ has $6$ edges, this cannot
happen; thus $C$ is closed, and therefore an equiangular hexagon.

  Finally, the lower adjacent vertex to
$v_2$ must be a boundary vertex $b_1$, to avoid long chains
starting at $b_L'$ or $b_R'$. Thus this connected component of
$\cal N$ is a hexagonal cell `anchored' to the
boundary.\vspace{.3cm}
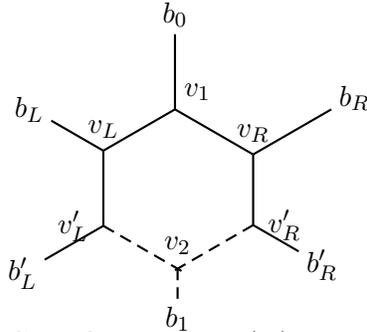
\begin{figure} \psset{unit=0.2cm}
\begin{center}
\begin{pspicture}(-8,-12)(8,12) \SpecialCoor 
\psline(0,0)(5;90) \psline(0,0)(5.5;210) \psline(0,0)(6;330)
\pnode(6;330){A} \psline(6;330)([angle=30,nodesep=6]A)
\psline(A)([angle=270,nodesep=4.7]A) \pnode(5.5;210){B}
\psline(B)([angle=150,nodesep=4]B)
\psline(B)([angle=270,nodesep=5]B)
\pnode([angle=270,nodesep=5]B){C}
\psline(C)([angle=210,nodesep=4.5]C)
\pnode([angle=270,nodesep=4.7]A){D}
\psline(D)([angle=330,nodesep=3.5]D)
\psline[linestyle=dashed](C)([angle=330,nodesep=6]C)
\psline[linestyle=dashed]([angle=210,nodesep=6]D)(D)
\pnode([angle=330,nodesep=5.7]C){E}
\psline[linestyle=dashed](E)([angle=270,nodesep=2]E)

\uput[ur](0,0){$v_1$} \uput{5.5}[90](0,0){$b_0$}
\uput[90](A){$v_R$} \uput[90](B){$v_L$} \uput{1}[180](C){$v_L'$}
\uput{1}[90](E){$v_2$}\uput{1}[0](D){$v_R'$}
\uput{0.5}[150]([angle=150,nodesep=4]B){$b_L$}\uput{0.5}[210]([angle=210,nodesep=4.5]C){$b_L'$}
\uput{0.5}[270]([angle=270,nodesep=2]E){$b_1$}\uput{0.5}[30]([angle=30,nodesep=6]A){$b_R$}
\uput{0.5}[330]([angle=330,nodesep=3.5]D){$b_R'$}

\end{pspicture}
 
\caption{Classification, part(iii)-hexagonal cell}\end{center}
\end{figure}

\emph{2.3 Connectedness of the network.}

Suppose $\cal N$ has two non-intersecting connected components
$\hat{\cal N}_1$ and $\hat{\cal N}_2$. In particular, the boundary
vertices of $\hat{\cal N}_1$ and $\hat{\cal N}_2$ are
`non-interlacing', that is, $\partial D$ is partitioned into two
arcs (disjoint except for their endpoints), each containing all
the boundary vertices of $\hat{\cal N}_1$ or $\hat{\cal N}_2$.
Each of $\hat{\cal N}_1$ and $\hat{\cal N}_2$ is of one of the
four types listed above (including critical chords). For any of
the four types, one sees directly that there is a connected arc of
$\cal C$ starting and ending at boundary vertices $B,\bar{B}$ of
that connected sub-graph, with a total turning angle from $B$ to
$\bar B$ (for the unit tangent to $\cal C$) of at least $\pi$
($\pi$ for a critical chord, $4\pi/3$ for a triode, $5\pi/3$ for a
double triode or a hexagonal cell.) For a strictly convex curve
$\cal C$, it is not possible for two such arcs to be disjoint.
Hence there is only one connected component.

We summarize the conclusion in the following proposition (Fig. 8):
\vspace{.2cm}

\textbf{Proposition 2.1.} \emph{Let $\cal N$ be a Steiner network
in a strictly convex domain in $\mathbb R^2$, with `free'
(Neumann) boundary conditions. Then $\cal N$ is connected, and is
one of: (i) a critical chord; (ii) a triode; (iii) a double
triode; (iv) a closed equiangular hexagon, anchored to the
boundary.}

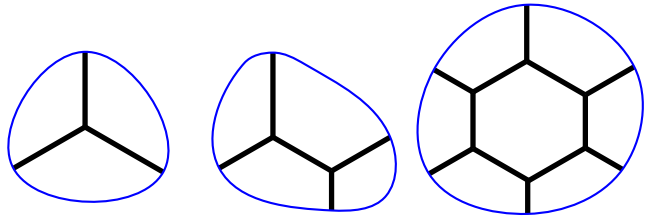
\begin{figure}
\begin{center}

\psset{unit=0.2cm} \SpecialCoor
\begin{pspicture}(14.5,4)(19.5,9)
 \psline[linewidth=2pt](0,0)(5;90) \psline[linewidth=2pt](0,0)(5.5;210)
\psline[linewidth=2pt](0,0)(6;330)
\psccurve[linecolor=blue](5;90)(5.5;210)(6;330)
\end{pspicture}

\psset{unit=0.15cm}
\begin{pspicture}(3,1)(9,6)
 \psline[linewidth=2pt](0,0)(7.5;90) \psline[linewidth=2pt](0,0)(5.5;210)
\psline[linewidth=2pt](0,0)(6;330) \pnode(6;330){V}
\psline[linewidth=2pt](6;330)([angle=30,nodesep=6]V)
\psline[linewidth=2pt](V)([angle=270,nodesep=3.5]V)
\pnode(7.5;90){A}\pnode(5.5;210){B}\pnode([angle=270,nodesep=3.5]V){C}\pnode([angle=30,nodesep=6]V){D}
\pnode([angle=200,nodesep=2.7]A){E}\pnode([angle=340,nodesep=3]A){F}\pnode([angle=15,nodesep=5]C){G}
\psccurve[linecolor=blue](A)(E)(B)(C)(G)(D)(F)
\end{pspicture}

\begin{pspicture}(-19,-11)(-14,-6)
\psset{unit=0.15cm}
\psline[linewidth=2pt](0,0)(5;90)
\psline[linewidth=2pt](0,0)(5.5;210)
\psline[linewidth=2pt](0,0)(6;330) \pnode(6;330){A}
\psline[linewidth=2pt](6;330)([angle=30,nodesep=5]A)
\psline[linewidth=2pt](A)([angle=270,nodesep=4.7]A)
\pnode(5.5;210){B}
\psline[linewidth=2pt](B)([angle=150,nodesep=4]B)
\psline[linewidth=2pt](B)([angle=270,nodesep=5]B)
\pnode([angle=270,nodesep=5]B){C}
\psline[linewidth=2pt](C)([angle=210,nodesep=4.5]C)
\pnode([angle=270,nodesep=4.7]A){D}
\psline[linewidth=2pt](D)([angle=330,nodesep=3.8]D)
\psline[linewidth=2pt](C)([angle=330,nodesep=5.6]C)
\psline[linewidth=2pt]([angle=210,nodesep=6]D)(D)
\pnode([angle=330,nodesep=5.6]C){E}
\psline[linewidth=2pt](E)([angle=270,nodesep=3]E)
\pnode(5;90){H1}\pnode([angle=150,nodesep=4]B){H2}\pnode([angle=210,nodesep=4.5]C){H3}
\pnode([angle=270,nodesep=3]E){H4}\pnode([angle=330,nodesep=3.8]D){H5}\pnode([angle=30,nodesep=5]A){H6}
\psccurve[linecolor=blue](H1)(H2)(H3)(H4)(H5)(H6)
\end{pspicture}

 
\caption{The three possible Steiner networks in a strictly convex
domain }\end{center}
\end{figure}


\newpage

\textbf{3. Existence of triodes.}

 \emph{ 3.1 Preliminary remarks}.
Given a strictly convex curve $\cal C$ and a direction $\theta \in
S$, let $\cal T(\theta)$ be the circumscribed equilateral triangle
with one side parallel to $T(\theta)$, touching $\cal C$ at the
points $B_1=B(\theta)$, $B_2=B(\theta+2\pi/3),
B_3=B(\theta+4\pi/3)$; for definiteness we take $\theta$ in
$[0,2\pi/3)$. The inner normals $n_1, n_2, n_3$ at the points of
contact intersect pairwise:
$$n_1\cap n_2=\{P_1\},\quad n_2\cap n_3=\{P_2\},\quad n_3\cap n_1=\{P_3\}$$
(Fig.9). $\theta$ defines a triode configuration exactly when
$P_1=P_2=P_3=P$ \emph{and} $P$ is in the interior of $D$. To
express this analytically, define functions $u(\theta)=u_1,
v(\theta)=v_1$ via:
$$P_1=B_1+\frac 2{\sqrt{3}}u_1N_1=B_2+\frac 2{\sqrt{3}}v_2N_2.$$
(Here $v_2=v(\theta+2\pi/3)$.) Taking inner products with $T_1$
and $T_2$ one finds:
$$u_1=\langle B_2-B_1, T_2\rangle,\quad v_2=\langle
B_2-B_1,T_1\rangle.$$ This gives the explicit definitions of the
`forward and backward triangle functions' $u(\theta), v(\theta)$:
$$u(\theta)=\langle
B(\theta+2\pi/3)-B(\theta),T(\theta+2\pi/3)\rangle,$$
$$v(\theta)=\langle
B(\theta)-B(\theta+4\pi/3),T(\theta+4\pi/3)\rangle.$$ Using
$T_1+T_2+T_3=0$, we see that $u-v$ is $2\pi/3$-periodic:
$$u(\theta)-v(\theta)=\langle B(\theta),T(\theta)\rangle+\langle
B(\theta+2\pi/3),T(\theta+2\pi/3)\rangle+\langle
B(\theta+4\pi/3),T(\theta+4\pi/3)\rangle.$$ Recalling (1.5) in
section 1, we see that $u-v=\omega$, the `derived height':
$$u(\theta)-v(\theta)=\omega(\theta)=h'(\theta).$$
This immediately implies the following. \vspace{.2cm}

\textbf{Proposition 3.1.} \emph{The three inner normals at the
points of contact of a circumscribed equilateral triangle $\cal
T(\theta)$ intersect at a single point exactly when $\theta\in S$
is a critical point of the height function; this happens for at
least two geometrically distinct configurations. A `critical
configuration' ($\omega(\theta)=0$) defines a triode iff
$u(\theta), u(\theta+2\pi/3), u(\theta+4\pi/3)$ are all positive.}

\begin{figure} \psset{unit=0.35cm}
\begin{center}

 \SpecialCoor
\begin{pspicture}(-5,-2)(15,18)
\pnode(0,0){Q3}  \pnode([angle=40,nodesep=20]Q3){Q1}
\pnode([angle=100,nodesep=20]Q3){Q2} \pspolygon(Q3)(Q1)(Q2)
\pnode([angle=40,nodesep=11]Q3){B1}
\pnode([angle=100,nodesep=11]Q3){B3}
\pnode([angle=160,nodesep=7]Q1){B2}
\psline[linecolor=green](B1)([angle=130,nodesep=10.5]B1)
\psline[linecolor=green](B2)([angle=250,nodesep=6.4]B2)
\psline[linecolor=green](B3)([angle=10,nodesep=9.75]B3)
\pnode([angle=210,nodesep=4]B1){B1a}
\pnode([angle=50,nodesep=4]B1){B1b}
\pnode([angle=330,nodesep=3]B2){B2a}
\pnode([angle=170,nodesep=3]B2){B2b}
\pnode([angle=90,nodesep=2]B3){B3a}
\pnode([angle=290,nodesep=2]B3){B3b}
\psccurve[linecolor=blue](B1)(B1b)(B2a)(B2)(B2b)(B3a)(B3)(B3b)(B1a)
\uput{0.5}[270](Q3){$Q_3$} \uput{0.5}[0](Q1){$Q_1$}
\uput{0.5}[180](Q2){$Q_2$} \uput{0.5}[0](B1){$B_1=B(\theta)$}
\uput{0.5}[45](B2){$B_2$} \uput{0.5}[200](B3){$B_3$}
\uput{0.4}[250]([angle=250,nodesep=6.4]B2){$P_1$}
\uput{0.5}[10]([angle=10,nodesep=9.75]B3){$P_2$}
\uput{0.5}[130]([angle=120,nodesep=5.5]B1){$P_3$}
\uput{0.2}[130]([angle=130,nodesep=10.5]B1){$B(\theta^*)$}
\pnode([angle=220,nodesep=1.5]B1){B1u}
\pnode([angle=220,nodesep=3]B1){B1v}
\psline[linewidth=0.5pt,arrows=|<->|](B1u)([angle=130,nodesep=3]B1u)
\psline[linewidth=0.5pt,arrows=|<->|](B1v)([angle=130,nodesep=6.4]B1v)
\uput{2}[140](B1u){$u$} \uput{5}[135](B1v){$v$}
\end{pspicture}
 
\caption{The forward and backward triangle functions }\end{center}
\end{figure}
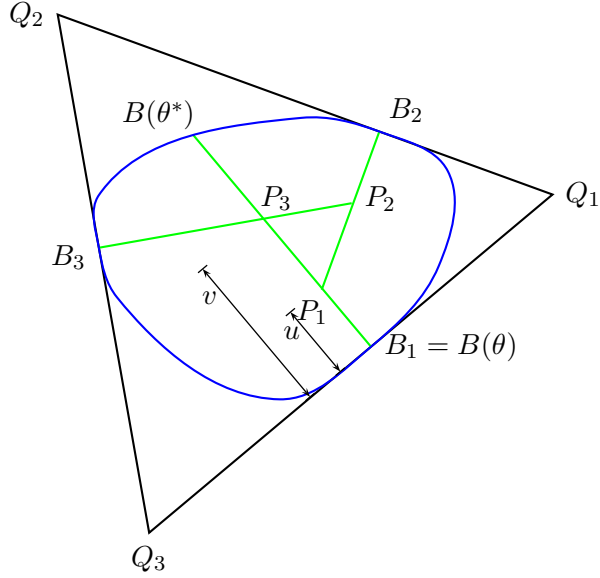

It is useful to express $u(\theta)$ in terms of the support
function $p$. From $\langle B_2,T_2\rangle=p'(\theta+2\pi/3)$ and:
$$\langle B_1,T_2\rangle=-p(\theta)\langle N_1,T_2\rangle
+p'(\theta)\langle T_1,T_2\rangle=-\frac{\sqrt{3}}2p(\theta)-
\frac 12 p'(\theta),$$ we have:
$$u(\theta)=p'(\theta+2\pi/3)+\frac 12
p'(\theta)+\frac{\sqrt{3}}2p(\theta).\qquad(3.1)$$ A similar
calculation yields:
$$v(\theta)=-p'(\theta+4\pi/3)-\frac 12
p'(\theta)+\frac{\sqrt{3}}2 p(\theta).$$ \vspace{.2cm}

\emph{Definition 3.1.} It is useful to observe that given a
(sufficiently differentiable) $2\pi$-periodic function $p$  one
can always define the `radius of curvature' $r=p''+p$, the
`derived height':
$$\omega(\theta)=p'(\theta)+p'(\theta+2\pi/3)+p'(\theta+4\pi/3)$$
and the `triangle function' $u$ (by (3.1)), even when $p$ is not
the support function of a strictly convex curve. \vspace{.2cm}

\emph{Remark 3.1.} Let $L(\theta)=B(\theta+2\pi/3)-B(\theta)$.
Then $\langle L(\theta),N(\theta)\rangle >0$ for all $\theta$ (by
strict convexity). Observing that $\sqrt{3}N_1=T_2-T_3$, and using
$u(\theta)=\langle L(\theta),T_1\rangle$,
$v(\theta+2\pi/3)=\langle L(\theta), T_1\rangle$, we have:
$$\begin{array}{ll}
u(\theta)>\langle L(\theta),T_3\rangle &=-\langle
L(\theta),T_1\rangle -\langle L(\theta),T_2\rangle\\
&=-v(\theta+2\pi/3)-u(\theta),
\end{array}$$
or $2u(\theta)+v(\theta+2\pi/3)>0$ for all $\theta$. By $2\pi/3$
periodicity of $u-v$, this is equivalent to:
$$u(\theta)+v(\theta)+u(\theta+2\pi/3)>0\quad \forall \theta.$$
A similar argument yields the equivalent inequalities:
$$2v(\theta)+u(\theta+4\pi/3)>0,\quad
u(\theta)+v(\theta)+v(\theta+2\pi/3)>0 \quad \forall \theta.$$
This shows that, for any $\theta$, at most two of $u_1,u_2,u_3$
can be negative; and at a zero of $\omega$, at most \emph{one} of
$u_1, u_2, u_3$ may be negative (for a triode, all must be
positive.) \vspace{.3cm}

\emph{3.1 Existence of triodes.} It is natural to look for
conditions of `curvature pinching' type that guarantee $u>0$
everywhere. \vspace{.2cm}

For each $\theta$, define $\theta_*\in(\theta,\theta+2\pi)$ by
requiring $B(\theta_*)$ to be the (other) intersection of the
normal through $B(\theta)$ with $\cal C$ (Fig.9). Set
$d(\theta)=\langle B(\theta_*)-B(\theta), N(\theta)\rangle$. All
intersections of inner normals (at contact points of a
circumscribed triangle $\cal T(\theta)$) will be internal if we
require, for all $\theta$:
$$0<u_1<\frac{\sqrt{3}}2d_1,0<u_2<\frac{\sqrt{3}}2d_2,0<u_3<\frac{\sqrt{3}}2d_3.\qquad (3.2)$$
For $C^2$ curves, we have:
$$d(\theta)=\int_0^{\theta_*-\theta}r(\tau+\theta)\sin \tau
d\tau$$ and:
$$u(\theta)=\int_{0}^{2\pi/3}r(\tau+\theta)\cos(2\pi/3-\tau)dt,$$
so we want:
$$0<-\frac 12\int_0^{2\pi/3} r(\tau+\theta)\cos \tau d\tau
+\frac{\sqrt{3}}2\int_0^{2\pi/3}r(\tau+\theta)\sin \tau d\tau <
\frac{\sqrt{3}}2\int_0^{\theta_*-\theta}r(\tau+\theta)\sin \tau
d\tau.$$ This suggests the condition:
$$\theta+2\pi/3<\theta_*<\theta+4\pi/3\quad \forall \theta.\qquad (3.3)$$
(Note that, for the circle, $\theta_*=\theta + \pi$  $\forall
\theta$.) This condition implies, in particular:
$$\theta+4\pi/3<(\theta+2\pi/3)_*<\theta+2\pi,\quad
\theta+2\pi<(\theta+4\pi/3)_*<\theta+8\pi/3\qquad \forall
\theta.$$ It is clear geometrically that:
$$u_1<0\Rightarrow (\theta+2\pi/3)_*>\theta+2\pi$$
and
$$u_1>(\sqrt{3}/2)d(\theta)\Rightarrow(\theta+2\pi/3)_*
<\theta_*<\theta+4\pi/3;$$ thus, (3.3) indeed implies (3.2).

It is easy to express (3.3) as a condition on the radius of
curvature. Since $\theta_*$ is characterized by:
$$\int_\theta^{\theta_*}r(\tau)\langle T(\tau),T(\theta)\rangle
d\tau=\int_{\theta}^{\theta_*}\cos(\tau-\theta)r(\tau)d\tau=0,$$
(2.3) is equivalent to the two conditions:
$$\int_0^{2\pi/3}r(\tau+\theta)\cos \tau d\tau>0,\quad
\int_0^{4\pi/3}r(\tau+\theta)\cos \tau d\tau<0.\qquad (3.4)$$
Assume $0<r_{min}\leq r(\theta)\leq r_{max}$ in $[0,2\pi]$. The
first inequality in (3.4) is equivalent to:
$$\int_0^{\pi/2}r(\tau+\theta)\cos \tau d\tau >
\int_{\pi/2}^{2\pi/3}r(\tau+\theta)|\cos \tau |d\tau,$$ which
would follow from $r_{min}>(1-\sqrt{3}/2)r_{max}.$ Similarly, the
second inequality in (3.4) is equivalent to:
$$\int_0^{\pi/2}r(\tau+\theta)d\tau<\int_{\pi/2}^{4\pi/3}r(\tau+\theta)|\cos
\tau| d\tau,$$ which would follow from:
$r_{max}<(1+\sqrt{3}/2)r_{min}$. This is the more restrictive
inequality, and gives the sufficient condition for existence of
triodes recorded in the following proposition.\vspace{.2cm}

\textbf{Proposition 3.2.} \emph{Let $p\in {\cal P}_{smooth}$
satisfy $r_{max}<(1+\sqrt{3}/2)r_{min}$. Then the three inner
normals at the contact points of a circumscribed equilateral
triangle always intersect inside $\cal C$. Thus $\cal C$ supports
at least two geometrically distinct triodes
(stably).}\vspace{.2cm}

\emph{Example 2.1.} The constant-height biangle has $u<0$ on the
`regular intervals' $(-\pi/6,\pi/6)$ and $(5\pi/6, 7\pi/6)$
(Fig.2). Indeed its `forward triangle function' $u_{bi}(\theta)$
is even and $\pi$-periodic, given in $[-\pi/6, 5\pi/6]$ by:
$$u_{bi}(\theta)=\left \{ \begin{array}{ll}&\frac{\sqrt{3}}2-\cos
\theta,\quad -\frac{\pi}6\leq \theta \leq \frac{\pi}6,\\&\frac
12(\sqrt{3}\sin \theta-\cos \theta),\quad \frac{\pi}6\leq \theta
\leq \frac{5\pi}6.\end{array}\right .$$ $u$ increases from the
minimum $\frac{\sqrt{3}}2-1$ at $\theta=0$ to the maximum
$\sqrt{3}/2$ at $\theta=\pi/2$, and is positive exactly on the
`singular intervals' $(\frac{\pi}6,\frac{5\pi}6)$ and
$(\frac{7\pi}6,\frac{11\pi}6)$, where $B(\theta)$ maps to one of
the endpoints of the maximal diameter. Thus the biangle does not
support triodes. We can remedy this by adding $c>0$ to the support
function, producing $p_c\in {\cal P}_{pw}$ (without corners) with
`triangle function' $u_c(\theta)=u_{bi}(\theta)+c$. For
$0<c<1-\frac{\sqrt{3}}2$, we have $u_c<0$ on the interval
$(-\delta,\delta)$, for some $0<\delta <\frac{\pi}6$.
($u_c(\delta)=0$). (Geometrically, the biangle and each of its
outer parallel curves have the \emph{same} locus of normal
intersections, shown in Fig.2) But if $\theta \in
(\delta,\frac{\pi}6)$, we have $\theta+\frac{2\pi}3\in
(\frac{2\pi}3+\delta,\frac{5\pi}6)$ and
$\theta+\frac{4\pi}3\in(\frac{4\pi}3+\delta,\frac{3\pi}2)$, and
$u_c>u_{bi}>0$ on both intervals; so $u_1,u_2,u_3$ are all
positive for these values of $\theta$, and since $\omega\equiv 0$
each value of $\theta$ in $(\delta, \frac{\pi}6)$ corresponds to a
triode. More generally, we have the following
proposition.\vspace{.2cm}

\textbf{Proposition 3.3.} \emph{Any $p\in {\cal P}_{pw}$ of
constant height supports triodes.}

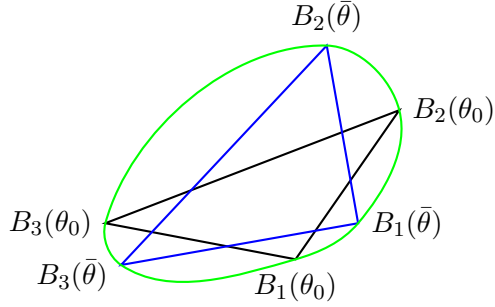
\begin{figure} \psset{unit=0.8cm}
\begin{center}
\SpecialCoor
\begin{pspicture}(-3,-1)(0,2)
\pnode(0,0){B1}\pnode(3;100){B2}\pnode(4;190){B3}\pnode(1.2;210){B10}
\pnode(2;70){B20}\pnode(4.2;180){B30} \pspolygon (B10)(B20)(B30)
\pspolygon[linecolor=blue](B1)(B2)(B3)
\psccurve[linecolor=green](B10)(B1)(B20)(B2)(B30)(B3)
\uput[270](B10){$B_1(\theta_0)$}\uput[0](B1){$B_1(\bar{\theta})$}\uput[0](B20){$B_2(\theta_0)$}
\uput[90](B2){$B_2(\bar{\theta})$}\uput[180](B30){$B_3(\theta_0)$}\uput[200](B3){$B_3(\bar{\theta})$}
\end{pspicture}
 
\caption{Proof of proposition 3.3}\end{center}
\end{figure}

\emph{Proof.} (Fig.10.) Since $\cal C$ is a $C^1$ strictly convex
curve, we may use a topological argument. For a given $\theta_0$,
consider the triangle of tangency points
$\Delta(\theta_0)=\{B_1(\theta_0),B_2(\theta_0),B_3(\theta_0)\}\subset
\bar{D};$ if all interior angles of $\Delta(\theta_0)$ are less
than $2\pi/3$, the normals to $\cal C$ at these points intersect
at a single point in $D$, and we are done. Assume, instead, that
the interior angle of $\Delta(\theta_0)$ at $B_1(\theta_0)$ is
greater than $2\pi/3$, while the angle at $B_2(\theta_0)$ is
smaller than $\pi/2$. By continuity, moving $B_1$ towards $B_2$ we
find $\bar{\theta}$ so that the angle of $\Delta(\bar{\theta})$ at
$B_1(\bar{\theta})$ is exactly $\pi/2$, so that \emph{all}
interior angles of $\Delta(\bar{\theta})$ are smaller than
$2\pi/3$. Hence the normals to $\cal C$ at the $B_i(\bar{\theta})$
intersect internally, and since $\omega \equiv 0$ we see
$\bar{\theta}$ defines a triode.\vspace{.2cm}

\emph{Example 3.2.} For the Reuleaux triangle, with support
function $p\in {\cal P}_{corner}$ given by (1.3), one checks
easily that $u\geq 0$ and vanishes only at $\theta=\pi/2, 3\pi/2$.
On the other hand, $\omega$ vanishes exactly at $0$ and $\pi/3$
(in $[0,2\pi/3)$), and both are transversal zeros of $\omega$.
$\theta=0, 2\pi/3, 4\pi/3$ all map to regular points of $\cal C$,
and hence define a triode, while $\theta=\pi/3,\pi,5\pi/3$ map to
the three vertices of the Reuleaux triangle, and again define a
triode (if we relax the definition slightly). For the `smoothed
Reuleaux triangles' defined by $p_{\epsilon} $ in (1.4), the
corresponding $\omega_{\epsilon}$ vanishes at the same points, and
$u_{\epsilon}$ is everywhere positive; so these curves also
support two triodes. Note that the ratio $r_{max}/r_{min}$ can be
made arbitrarily large, by taking $\epsilon>0$ small.
\vspace{.2cm}

We end this subsection by recording four simple observations
regarding existence for triodes: (i) any $\cal C$ corresponding to
$p\in {\cal P}_{pw}$ with an axis of reflection symmetry supports
triodes. (ii) If $\cal C$ supports a triode, it extends to a triode
in any exterior parallel curve; (iii) For \emph{any} strictly convex
$\cal C$, all exterior parallel curves ${\cal C}_c$ for $c$
sufficiently large will support triodes (since adding a sufficiently
large positive constant to $p$ makes $u$ everywhere positive.) (iv)
Curves of constant width ($w\equiv 0$) always support triodes (see
Remark 5.5 in section 5.)\vspace{.2cm}

\emph{3.3 Strictly convex curves without triodes.}

We use the nonexistence criterion: if $u<0$ at every zero of
$\omega$ in a `fundamental domain' (i.e., an interval of length
$2\pi/3)$, then $\cal C$ supports no triodes. We can achieve this
by starting with a support function for a strictly convex domain
of `constant height' satisfying $u<0$ in some interval, and taking
convex combinations with another `support function' for which the
zeros of $\omega$ have the desired property. For the second one,
we don't even need $r>0$. \vspace{.2cm}

\textbf{Proposition 3.4.} \emph{Let $p_0\in {\cal P}_{smooth}$
satisfy $\omega_0\equiv 0$ (`constant height') and $u_0<0$ in some
interval $I\subset (-\pi,\pi)$. Let $\tilde{p}$ be $2\pi$-periodic
and $C^2$, with the property that all zeros of the corresponding
`derived height' $\tilde{\omega}$ in a fundamental domain are
contained in an interval $J\subset \subset I$. Then for all
$\lambda \in (0,1)$ sufficiently close to $1$, the function:
$$p_{\lambda}(\theta)=(1-\lambda)\tilde{p}(\theta)+\lambda
p_0(\theta)$$ is in ${\cal P}_{smooth}$, and is the support
function of a $C^2$ strictly convex curve which does not support
triodes. If the zeros of $\tilde{\omega}$ are transversal, this is
true stably.}

\vspace{.2cm}\emph{Proof.} Since
$r_{\lambda}=(1-\lambda)\tilde{r}+\lambda r_0$, we clearly have
$p_{\lambda}\in {\cal P}_{smooth}$ (i.e., $p>0$ and $r>0$) for
$\lambda$ sufficiently close to 1. Also,
$u=(1-\lambda)\tilde{u}+\lambda u_0$, so for $\lambda$ close to 1
we have $u_{\lambda}<0$ in $I$; while
$\omega_{\lambda}=(1-\lambda)\tilde{\omega}$, so for \emph{all}
$\lambda\in (0,1)$ we have $\omega_{\lambda}<0$ in $J\subset I$.
So for $\lambda$ close to 1, all zeros of $\omega_{\lambda}$ in a
fundamental domain are contained in an interval where
$u_{\lambda}<0$. Hence $p_{\lambda}$ cannot support a triode.
\vspace{.2cm}

It is not hard to construct functions $\tilde{p}$ satisfying the
conditions of the proposition. For example, for any $0<m<1$, the
function:
$$q(\theta)=\frac{\cos(3\theta)}{1-m\sin(3\theta)}$$
is $2\pi/3$-periodic, with critical points exactly where
$sin(3\theta)=m$. Thus the zeros of $q'$ in $[0,2\pi/3)$ are
exactly $\frac 13 \arcsin (m)$ and $\frac{\pi}3-\frac 13 \arcsin(
m)$ (where $\arcsin (m)\in (0,\pi/2)$). The distance between these
zeros is $\frac 23 (\frac{\pi}2-\arcsin (m))$, which can be made
arbitrarily small by taking $m$ close enough to 1. Thus, given
$I\subset (-\pi,\pi)$ (defined by $p_0$ as in the proposition), a
suitable translate of $q$ can be used as $\tilde p$. (Note that,
since $q$ is $2\pi/3$-periodic, the corresponding derived height
function is simply $3q'$.)\vspace{.2cm}

\textbf{Corollary 3.5.} \emph{Let $p_0\in {\cal P}_{smooth}$ be
any curve of constant height, such that $u_0<0$ somewhere. Then
arbitrarily $C^r$ close to $p_0$ (if $p_0\in C^r$, where $r\geq
2$), one finds: (i) $C^r$ strictly convex curves supporting
triodes (stably); (ii) $C^r$ strictly convex curves which do not
support any triodes (also stably).} \vspace{.2cm}

\emph{Proof.} Part (ii) follows directly from the proposition; for
stability of the no-triode property, we just need to remark that
the zeros of $q$ given above are transversal. Clearly, by taking
$\lambda$ sufficiently close to 1, $p_{\lambda}$ can be made
arbitrarily $C^r$-close to $p_0$. For (i), recall that by Prop.
(3.3), one may find $\bar{\theta}$ so that $u_1^0, u_2^0, u_3^0$
(computed for $p_0$) are all positive at $\bar{\theta}$. As
explained above, we find $\tilde{p}$ smooth, $2\pi/3$-periodic, so
that all the zeros of the corresponding $\tilde{\omega}$ in a
fundamental domain are found in a small neighborhood of
$\tilde{\theta}$, where $u_1^0, u_2^0, u_3^0$ are still positive.
Taking now:
$$p_{\lambda}=(1-\lambda)\tilde{p}+\lambda p_0,$$
we have, for $\lambda$ sufficiently close to 1: (i)
$p_{\lambda}\in {\cal P}_{smooth}$ and is as close to $p_0$ as
desired; (ii) the zeros of
$\omega_{\lambda}=(1-\lambda)\tilde{\omega}$ in a fundamental
domain are all found in a neighborhood of $\tilde{\theta}$ where
$u_1^{\lambda},u_2^{\lambda},u_3^{\lambda}$ are all positive. Thus
$p_{\lambda}$ supports triodes.\vspace{.3cm}

\emph{Example 3.3.} An explicit example of $p_0\in {\cal
P}_{smooth}$ of constant height with $u<0$ somewhere is the
`smoothed biangle' of (1.7):
$$p_{\epsilon}(\theta)=\frac{\sqrt{3}}\pi(\frac 34-\frac
13\cos(2\theta)-\frac 1{30}\cos(4\theta))+\epsilon$$ ($\epsilon
>0$ arbitrary). The `triangle function' $u_{\epsilon}$ computed
from $p_{\epsilon}$ is:
$$u_{\epsilon}(\theta)=\frac 3{2\pi}(-\cos(2\theta)+\frac
1{10}\cos(4\theta)+\frac 34)+\frac{\sqrt{3}}2 \epsilon.$$ For
$\epsilon=0$, this is negative in $(-x_0,x_0)$ and vanishes at
$$x_0=(1/2)\arccos(5/2-\sqrt{3})=0.3476\ldots$$  Since
$u_{\epsilon}=u_0+\frac{\sqrt{3}}2\epsilon$ and $u_0(0.34)\sim
-\frac{0.01}\pi$, we find that for $\epsilon=0.02/(\pi\sqrt{3})$,
we have $u_{\epsilon}<0$ in $I=(-0.34,0.34)$.

Turning to parameters in $\tilde {p}$, we let
$m=\sin(3\pi/4)=\sqrt{2}/2$; then the critical points of $q(x)$ in
$[0,2\pi/3)$ are $\pi/4=0.785\ldots$ and $\pi/12=0.262\ldots$;
letting $\tilde{p}(\theta)=q(\theta+\pi/6)$, $\tilde{\omega}$ has
in the fundamental domain $[-\pi/6,\pi/2)$ only the zeros
$-\pi/12=-0.262\ldots$ and $\pi/12=0.262\ldots$, both in the
interval $I$. \vspace{.2cm}

\emph{Remark 3.2.} One checks numerically that the Fourier series
truncated at $n=9$ of
$\tilde{p}(\theta)=q(\theta+\pi/6)=-\sin(3\theta)/(1-(\sqrt{2}/2)\cos(3\theta))$
already has the desired property. Writing $p_{\lambda}$ in the
form (for $p_{\epsilon}$ and $\epsilon$ given above:)
$$p_{\lambda}(\theta)=\frac{\lambda}{1+\lambda}p_{\epsilon}(\theta)+
\frac 1{(1+\lambda)}\tilde{p}(\theta),$$ one finds that for
$\lambda=5,000$ we already have $p_{\lambda}>0$ and
$r_{\lambda}>0$ for all $\theta$. The resulting convex curve is
visually indistinguishable from the biangle.\vspace{.2cm}

\emph{Remark 3.3.} The construction depends essentially on the
existence of convex curves of constant height for which the
`triangle function' $u$ is negative somewhere. While this happens
for the constant height biangle, the other `curved regular polygons'
of constant height have $u>0$ everywhere, and any convex curve of
constant height may be uniformly approximated by such polygons. So
our construction of examples of convex curves without triodes
ultimately depends on the existence of this atypical example of a
convex constant-height curve.\vspace{.2cm}

\emph{3.4. Convergence of triode configurations.}\vspace{.2cm}

\emph{Definition 3.2.} $\theta \in S$ defines a \emph{boundary
triode} if $\omega(\theta)=0$ and $u(\theta)=0$; geometrically,
one of the edges of the triode collapsed to a boundary point of
$D$.

As long as the limit convex curve is without corners, this is the
only kind of degeneration allowed under convergence in ${\cal
P}_{smooth}$ with the $C^2$ topology. \vspace{.2cm}

\textbf{Lemma 3.6.} If $p\in {\cal P}_{pw}$ with (piecewise
continuous) radius of curvature function $r(\theta)\geq r_0>0$,
then for each $\theta \in \mathbb R$:
$$||B(\theta+2\pi/3)-B(\theta)||\geq (3/2)\min \{r(\tau);\tau \in
[\theta,\theta+2\pi/3]\}\geq (3/2)r_0.$$

\emph{Proof.} Since:
$$B(\theta+2\pi/3)-B(\theta)=(\int_0^{2\pi/3}r(\tau+\theta)\cos
\tau d\tau)T(\theta)+(\int_0^{2\pi/3}r(\tau +\theta)\sin \tau
d\tau)N(\theta),$$ we have:
$$||B(\theta+2\pi/3)-B(\theta)||^2\geq (\int_0^{2\pi/3}r(\tau
+\theta)\sin \tau d\tau)^2\geq
(3/2)^2(\min_{[\theta,\theta+2\pi/3]}r)^2.$$ \vspace{.2cm}

\textbf{Proposition 3.7.} Let $p_i,p\in {\cal P}_{smooth}$,
$p_i\rightarrow p$ uniformly in $C^2(S)$. Let ${\mathbb
T}_i=(B_1^i,B_2^i,B_3^i,P^i)$ be a triode or boundary triode in
${\cal C}_i$, where $P^i\in \bar{D}_i$. Then (up to passing to
subsequences) ${\mathbb T}_i\rightarrow {\mathbb
T}=(B_1,B_2,B_3,P)$ (meaning $B_1^i\rightarrow B_1, P^i\rightarrow
P$, etc.)

\emph{Proof.} This is practically self-evident, since
$\omega_i(\theta_i)=0$ implies $\omega(\theta)=0$. From the lemma,
$||B^i_a-B^i_b||$ is uniformly bounded below for $a\neq b$ in
$\{1,2,3\}$, so the only `collapse' allowed is
$||B_a^i-P^i||\rightarrow 0$ for some $a$ (and then
$P^i\rightarrow P\in {\cal C}$), which gives a boundary triode in
the limit. \vspace{.2cm}

\emph{Remark 3.4.} If $p\in {\cal P}_{corner}$ and $\cal C$ has
corners, it is conceivable that $p_i\rightarrow p$ (say, uniformly
in $C^1$), but the $\mathbb T_i$ collapse to a segment. This
happens, for example, when triodes in the outer parallel curves of
a biangle converge to the diameter of the biangle. They cannot
collapse to a point, however, since the total length $L_i$ is
uniformly bounded below, as long as $p_i(\theta)\geq p_0>0$.


\newpage
\textbf{4. Double triode configurations.}\vspace{.2cm}

\emph{4.1 Basic construction.} We adopt for the remainder of the
paper the notation $\bar{\theta}=\theta+\pi$. Let $\cal C$ be
strictly convex. Given any direction $\theta \in S$, the
circumscribed triangles (given by their vertices) ${\cal
T}(\theta)=\{Q_1,Q_2,Q_3\}$ and ${\cal
T}(\bar{\theta})=\{\bar{Q_1},\bar{Q_2},\bar{Q_3}\}$ define three
parallelograms:
$${\cal P}(\theta)=\{Q_1,S_2,\bar{Q_1},\bar{S_2}\}={\cal
P}({\bar{\theta}})$$ ( Fig.11) and also (see Fig. 15):
$${\cal P}(\theta+2\pi/3)=\{Q_2,S_3,\bar{Q}_2,\bar{S}_3\},\quad
{\cal P}(\theta+4\pi/3)=\{Q_3,S_1,\bar{Q}_3,\bar{S}_1\}.$$
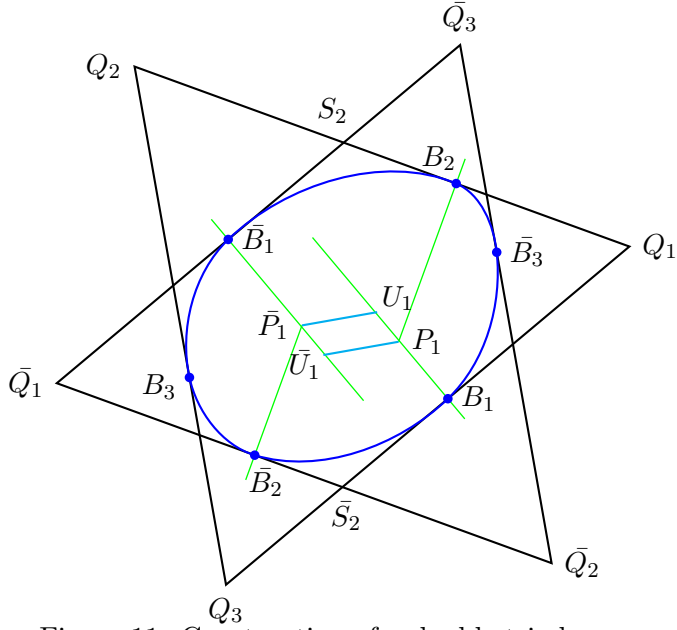
\begin{figure} \psset{unit=0.35cm}
\begin{center}
\SpecialCoor
\begin{pspicture}(-6,0)(12,18)
\pnode(0,0){Q3}  \pnode([angle=40,nodesep=20]Q3){Q1}
\pnode([angle=100,nodesep=20]Q3){Q2} \pspolygon(Q3)(Q1)(Q2)
\pnode([angle=40,nodesep=11]Q3){B1}
\pnode([angle=100,nodesep=8]Q3){B3}
\pnode([angle=160,nodesep=7]Q1){B2}
\psline[linewidth=0.5pt,linecolor=green]([angle=310,nodesep=1]B1)([angle=130,nodesep=8]B1)
\psline[linewidth=0.5pt,linecolor=green]([angle=70,nodesep=1]B2)([angle=250,nodesep=6.4]B2)
\pnode(10;130){Q1bar} \pnode([angle=340,nodesep=20]Q1bar){Q2bar}
\pnode([angle=40,nodesep=20]Q1bar){Q3bar}
\pspolygon(Q1bar)(Q2bar)(Q3bar)
\pnode([angle=340,nodesep=8]Q1bar){B2bar}
\pnode([angle=40,nodesep=8.5]Q1bar){B1bar}
\pnode([angle=100,nodesep=12]Q2bar){B3bar}
\psline[linewidth=0.5pt,linecolor=green]([angle=70,nodesep=5.25]B2bar)([angle=250,nodesep=1]B2bar)
\psline[linewidth=0.5pt,linecolor=green]([angle=310,nodesep=8]B1bar)([angle=130,nodesep=1]B1bar)
\pnode([angle=70,nodesep=5.25]B2bar){P1bar}
\pnode([angle=250,nodesep=6.4]B2){P1}
\psline[linecolor=cyan](P1bar)([angle=10,nodesep=2.9]P1bar)
\pnode([angle=10,nodesep=2.9]P1bar){U1}
\psline[linecolor=cyan](P1)([angle=190,nodesep=2.9]P1)
\pnode([angle=190,nodesep=2.9]P1){U1bar}
\psccurve[showpoints=true,linecolor=blue](B2bar)(B1)(B3bar)(B2)(B1bar)(B3)
\uput{0.5}[270](Q3){$Q_3$} \uput{0.5}[0](Q1){$Q_1$}
\uput{0.5}[180](Q2){$Q_2$}\uput{0.5}[0](B1){$B_1$}
\uput{0.5}[120](B2){$B_2$}\uput{0.5}[0](P1){$P_1$}\uput{0.2}[40](U1){$U_1$}
\uput{0.5}[200](B3){$B_3$} \uput{0.5}[180](Q1bar){$\bar{Q_1}$}
\uput{0.5}[0](Q2bar){$\bar{Q_2}$}
\uput{0.5}[90](Q3bar){$\bar{Q_3}$}
\uput{0.5}[0](B1bar){$\bar{B_1}$}
\uput{0.5}[290](B2bar){$\bar{B_2}$}
\uput{0.5}[0](B3bar){$\bar{B_3}$}
\uput{0.5}[180](P1bar){$\bar{P_1}$}
\uput{0.1}[210](U1bar){$\bar{U_1}$}
\uput{4.5}[30](Q3){$\bar{S_2}$} \uput{14}[45](Q1bar){$S_2$}

\end{pspicture}
 
\caption{Construction of a double triode}\end{center}
\end{figure}
For each $\theta\in S$, intersections of normals define a four-point
configuration
$(P(\theta),U(\theta),P(\bar{\theta}),U(\bar{\theta}))$, collapsing
to two points when $\theta$ defines a double triode. Here
$P(\theta)=n(\theta)\cap n(\theta+2\pi/3)$, as in section 3. To
define $U(\theta)$, we consider `projections'. Since there are only
three unoriented directions involved in any potential network
configuration (once $\theta$ is fixed), given any two non-parallel
normal lines $n$, $\tilde{n}$ along two of the directions, there is
a well-defined projection from $n$ to $\tilde{n}$ along the third
normal direction, denoted $pr_{n\rightarrow \tilde{n}}$. For a
double triode configuration, we consider the points:
$$U(\theta)=pr_{\bar{n}_2\rightarrow n_1}(P(\bar{\theta})),\quad
\bar{U}(\theta)=U(\bar{\theta})=pr_{n_2\rightarrow\bar{n}_1}(P(\theta))$$
(where $n_1=n(\theta),\bar{n}_2=n(\theta+2\pi/3+\pi)$, etc.) We
write this in simplified notation as follows:
$$U_1=pr_{\bar{2}1}(\bar{P}_1),\quad \bar{U}_1=pr_{2\bar{1}}(P_1).$$
We have a double triode when (i) $P_1=U_1$ (and then
$\bar{P}_1=\bar{U}_1$); and (ii) $P_1$and $\bar{P}_1$ are inside
$\cal C$. (A third necessary condition will be described soon).
Then the parallelogram $P_1\bar{U}_1\bar{P}_1U_1$ collapses to the
`bridge' of the double triode ${\mathbb X}=\{B(\theta),
B(\theta+2\pi/3),B(\bar{\theta}),B(\bar{\theta}+2\pi/3)\}$.

We let $u(\theta)=u_1$ have the same meaning as in section 3, and
define $\sigma(\theta)=\sigma_1,
\sigma(\bar{\theta})=\bar{\sigma}_1$:
$$P_1=B_1+\frac 2{\sqrt 3}u_1N_1, U_1=B_1+\frac 2{\sqrt
3}\sigma_1 N_1, \bar{P}_1=\bar{B}_1+\frac 2{\sqrt 3}
\bar{u_1}\bar{N_1}, \bar{U_1}=\bar{B}_1+\frac 2{\sqrt
3}\bar{\sigma}_1\bar{N_1}.$$

Since $U_1$ is the projection of $\bar{P_1}$ along the direction
$N_3$, we also have: $U_1=\bar{P_1}+\alpha N_3$ for some
$\alpha\in \mathbb R$; taking inner products with $T_3$, we find:
$$\sigma_1=\langle B_1-\bar{B_1},T_3\rangle-\bar{u_1}.$$
We are interested in the difference $\sigma-u$:
$$\begin{array}{rcl}
\sigma-u&=&\langle B_1-\bar{B}_1,T_3\rangle-\langle
\bar{B}_2-\bar{B}_1,\bar{T_2}\rangle-\langle B_2-B_1,T_2\rangle\\
&=&\langle \bar{B}_1-B_1,T_1\rangle+\langle
\bar{B}_2-B_2,T_2\rangle\\
&=&w(\theta)+w(\theta+2\pi/3),\end{array}$$ where the derived
width function $w$ was defined in section 1.  We record the
expression for $\sigma$:
$$\sigma(\theta)=u(\theta)+w(\theta)+w(\theta+2\pi/3)\quad (4.1)$$
Condition (i) above for existence of a double triode is
$u=\sigma$, or:
$$ y(\theta):=w(\theta)+w(\theta+2\pi/3)=0.$$
We can relate this to the perimeter $L_{par}$ of the circumscribed
parallelogram ${\cal P}(\theta)$:
$$\begin{array}{rcl}L_{par}(\theta)&=&2(\langle S_2-Q_1,T_2\rangle+\langle
\bar{Q}_1-S_2,\bar{T_1}\rangle)\\
&=&\frac 4{\sqrt 3}(b(\theta)+b(\theta+2\pi/3)).\end{array}$$
Since $b'(\theta)=-w(\theta)$, we conclude:
$$y(\theta)=0\Leftrightarrow \theta\mbox{ is a critical point of
the perimeter of }{\cal P}(\theta).$$ \vspace{.1cm}

It is useful to observe that, since $w$ is $\pi$-periodic, one can
recover $w$ from $y$ via:
$$y(\theta)+y(\theta+\pi/3)-y(\theta+2\pi/3)=2w(\theta),\quad (4.2)$$
with a similar relation between $L_{par}$ and $b$. Thus the class
`curves of constant width' is the same as `curves of constant
parallelogram perimeter'. \vspace{.2cm}

The third condition for a double triode arises from the fact that
two types of configurations are possible if $y(\theta)=0$.
(Fig.12).
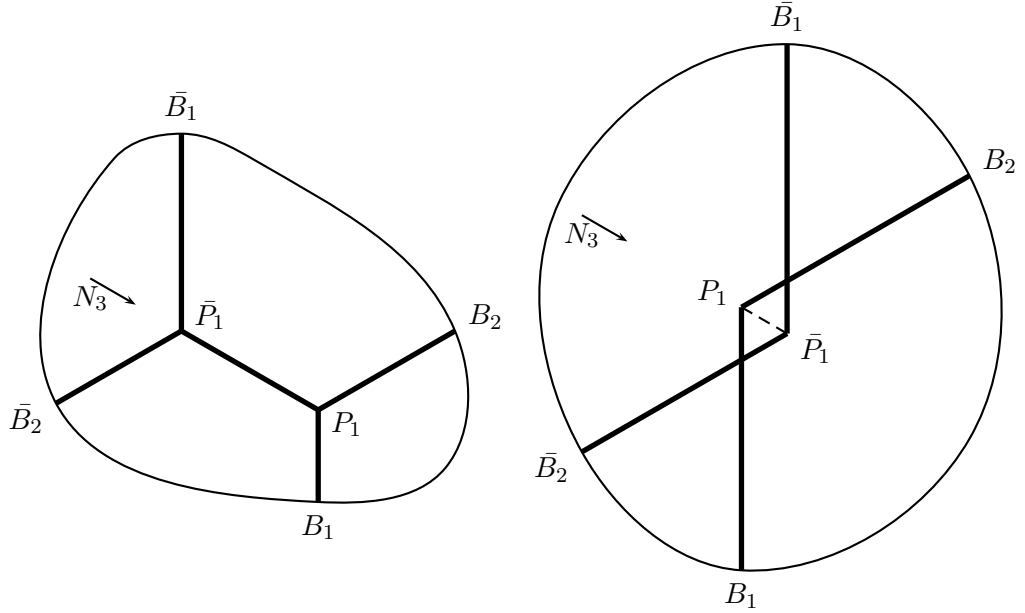
\begin{figure} \psset{unit=0.35cm}
\begin{center}
\SpecialCoor
\begin{pspicture}(5,-1)(15,9)
\psline[linewidth=2pt](0,0)(7.5;90)
\psline[linewidth=2pt](0,0)(5.5;210)
\psline[linewidth=2pt](0,0)(6;330) \pnode(6;330){V}
\psline[linewidth=2pt](6;330)([angle=30,nodesep=6]V)
\psline[linewidth=2pt](V)([angle=270,nodesep=3.5]V)
\psline{->}(4;150)(2;150) \uput{0.3}[270](4;150){$N_3$}
\pnode(7.5;90){A}\pnode(5.5;210){B}\pnode([angle=270,nodesep=3.5]V){C}\pnode([angle=30,nodesep=6]V){D}
\pnode([angle=200,nodesep=2.7]A){E}\pnode([angle=340,nodesep=3]A){F}\pnode([angle=15,nodesep=5]C){G}
\psccurve(A)(E)(B)(C)(G)(D)(F)
\uput[90](A){$\bar{B_1}$} \uput[210](B){$\bar{B_2}$}
\uput[270](C){$B_1$} \uput[30](D){$B_2$}
\uput[330](V){$P_1$}\uput[30](0,0){$\bar{P_1}$}
\end{pspicture}

\begin{pspicture}(-8,-1)(-18,-11)
\psline[linewidth=2pt](0,0)(11;90)
\psline[linewidth=2pt](0,0)(9;210)
\psline[linestyle=dashed](0,0)(2;150) \pnode(2;150){V}
\psline[linewidth=2pt](V)([angle=30,nodesep=10]V)
\psline[linewidth=2pt](V)([angle=270,nodesep=10]V)
\psline{->}(9;150)(7;150) \uput{0.3}[270](9;150){$N_3$}
\pnode(11;90){B1bar}\pnode(9;210){B2bar}\pnode([angle=30,nodesep=10]V){B2}
\pnode([angle=270,nodesep=10]V){B1}
\pnode([angle=150,nodesep=8]V){Q}
\pnode([angle=330,nodesep=10]V){R}
\psccurve(B2bar)(B1)(R)(B2)(B1bar)(Q)
\uput[330](0,0){$\bar{P_1}$} \uput[150](V){$P_1$}
\uput[90](B1bar){$\bar{B_1}$}\uput[210](B2bar){$\bar{B_2}$}
\uput[270](B1){$B_1$}\uput[30](B2){$B_2$}
\end{pspicture}

 
\caption{Double triode construction: good case, bad
case}\end{center}
\end{figure}

\emph{Case(I)} (good): $w(\theta)<0$ and $w(\theta+2\pi/3)>0$;

\emph{Case (II)} (bad): $w(\theta)>0$ and $w(\theta+2\pi/3)<0$.
\vspace{.2cm}

In fact, one easily computes that the (oriented) length of the
`bridge' of a critical configuration (which is positive for a
double triode) is: $\langle P_1-\bar{P_1},N_3\rangle=-(2/\sqrt
3)w(\theta)$, and the total length of the network is:
$$L=\frac 2{\sqrt{3}}(u_1+\bar{u}_1+v_2+\bar{v}_2-w).$$
Using the already computed expressions for the various terms, we
find:
$$L=b(\theta)+b(\theta+2\pi/3),$$
in words: the length of the double triode $\mathbb X(\theta)$ is
$\sqrt{3}/4$ times the perimeter of the parallelogram $\cal
P(\theta)$. \vspace{.2cm}

From (4.2) we see that a double triode corresponds to
$y(\theta)=0$ with $y(\theta+\pi/3)<y(\theta+2\pi/3)$. It follows
that 3-symmetric curves (for which $w$ and $y$ are both $\pi$- and
$2\pi/3$-periodic, hence $\pi/3$-periodic) cannot support
non-degenerate double triodes. The same holds for curves of
constant width ($w\equiv 0$).\vspace{.2cm}

To summarize: $\theta$ corresponds to a double triode
configuration if (i) $y(\theta):=w(\theta)+w(\theta+2\pi/3)=0$;
(ii) $u(\theta)>0$ and $u(\bar{\theta})>0$; (iii)$w(\theta)<0$.
Note that if $\theta$ is a transversal zero of $y$, any
sufficiently $C^2$ close curve will also support a double triode.
\vspace{.3cm}

\emph{4.2 Sufficient conditions for existence.} It is natural to
consider convex curves admitting an axis of reflection symmetry,
and then look for symmetric double triodes. Such an axis is always
a critical chord of $\cal C$- its endpoints correspond to critical
points of the width function $b$, which will be assumed to be
non-degenerate, hence a local max (`maximal chord') or local min
(`minimal chord') of $b$. Symmetry implies the existence of a
second critical chord orthogonal to the axis of symmetry. For our
first existence result we assume these are the only critical
chords of $\cal C$.\vspace{.2cm}

\textbf{Proposition 4.1.} \emph{Let $\cal C$ be a $C^1$ strictly
convex curve ($p\in {\cal P}_{pw}$) Assume $\cal C$ has a
non-degenerate maximal chord which is a line of reflection
symmetry and, except for the orthogonal minimal chord (also
assumed non-degenerate), no other critical chords. Then $\cal C$
stably supports a `double triode'.}

\emph{Examples} are  given by ellipses and outer parallel curves
of the constant-height biangle. \vspace{.2cm}
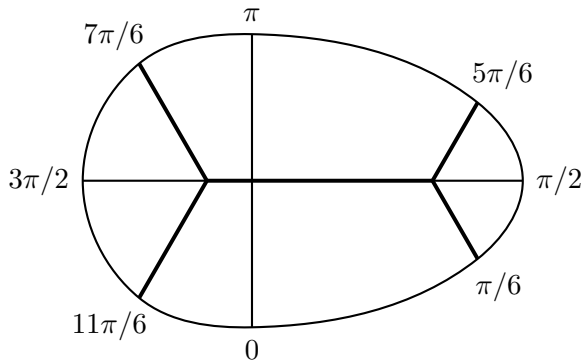
\begin{figure} \psset{unit=0.3cm}
\begin{center}
\SpecialCoor
\begin{pspicture}(-5,-7)(5,3)
\pnode(0,0){O}
\pnode(12;0){A}\pnode(6.5;90){B}\pnode(7.5;180){C}\pnode(6.5;270){D}
\pnode(8;0){P}\pnode(2;180){Pbar}
\pnode([angle=300,nodesep=4]P){B1}\pnode([angle=60,nodesep=4]P){B2}
\pnode([angle=120,nodesep=6]Pbar){B1bar}
\pnode([angle=240,nodesep=6]Pbar){B2bar}
\psline(B)(D)\psline(A)(C) \psline[linewidth=1.5pt](P)(Pbar)
\psline[linewidth=1.5pt](P)(B1) \psline[linewidth=1.5pt](P)(B2)
\psline[linewidth=1.5pt](Pbar)(B1bar)
\psline[linewidth=1.5pt](Pbar)(B2bar)
\psccurve(B2bar)(D)(B1)(A)(B2)(B)(B1bar)(C)
\uput[270](D){0}\uput[300](B1){$\pi/6$}\uput[0](A){$\pi/2$}\uput[60](B2){$5\pi/6$}
\uput[90](B){$\pi$}\uput[120](B1bar){$7\pi/6$}\uput[180](C){$3\pi/2$}
\uput[240](B2bar){$11\pi/6$}
\end{pspicture}
 
\caption{Existence of a double triode with symmetry
(Prop.4.1)}\end{center}
\end{figure}

\emph{Proof.} (Fig. 13) Position $\cal C$ so that the line of
symmetry has direction $T(0)$. Then the maximal chord is
$B(\frac{3\pi}2)B(\frac{\pi}2)$, the minimal chord $B(0)B(\pi)$.
$w=-b'$ is negative on $(0,\frac{\pi}2)\cup (\pi,\frac{3\pi}2)$, in
particular at $\theta_0=\pi/6$. It is easy to check that
$\theta_0=\pi/6$ satisfies the two other conditions for a double
triode. First, $p'(\pi-\theta)=-p'(\theta)$ (symmetry) implies
$w(\pi-\theta)=-w(\theta)$, and hence
$y(\frac{\pi}3-\theta)=-y(\theta)$, so $y(\pi/6)=0$. Second,
symmetry implies $B(5\pi/6)-B(\pi/6)=cN(0)$, for some $c>0$. Hence
$u(\frac{\pi}6)=c\cos (\frac{5\pi}6-\frac{\pi}2)>0$, and similarly
for $u(\frac{7\pi}6)$. Stability follows from the fact that $\pi/6$
is a transversal zero of $y$.\vspace{.2cm}

The next result assumes reflection symmetry only near the
endpoints of a maximal chord, and makes quantitative the intuition
that double triodes are easier to find near `sharp tips'. A
motivating example is a strictly convex domain with two corners,
so that at each corner no wedge with aperture greater than or
equal to $\pi/6$ about the chord joining the corners fits inside
$D$; it is easy to see that any outer parallel curve ${\cal
C}_{\epsilon}$ supports a double triode.

 To set up the notation,
let $\bar{K}K$ be a non-degenerate maximal chord with direction
$(1,0)$, so $K=B(\pi/2), \bar{K}=B(3\pi/2)$. We seek a double
triode with `bridge' along $\bar{K}K$  (that is, defined by
$\theta=\pi/6$), and \emph{assume} $y(\pi/6)=0$ (this would follow
from global reflection symmetry, as seen above.) Set:
$$r_K=\sup\{r(\theta); \theta \in [\pi/6, 5\pi/6]\}; \quad
r_{\bar{K}}=\sup \{r(\theta); \theta \in [7\pi/6, 11\pi/6]\}.$$
\vspace{.2cm}

\textbf{Proposition 4.2.} \emph{Let $\cal C$ be strictly convex
(say, defined by $p\in {\cal P}_{pw}$), with a non-degenerate
maximal chord $\bar{K}K$ of length $L_{chord}$. Assume (i)
$p(\pi-\theta)=p(\theta)$ for $\theta \in
[\frac{\pi}6,\frac{5\pi}6]\cup[\frac{7\pi}6,\frac{11\pi}6]$
(reflection symmetry near $K,\bar{K}$); (ii) $\theta_0=\pi/6$ is a
transversal zero of $y$; (iii) with $c_0=(4/3)(1+4\pi^2/3)^{1/2}$,
assume: $r_K+r_{\bar{K}}<\frac 1{c_0}L_{chord}$ ($c_0^{-1}\sim
0.2$.) Then $\cal C$ stably supports a double triode.}
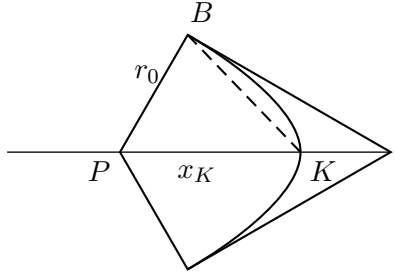
\begin{figure} \psset{unit=0.6cm}
\begin{center}
\SpecialCoor
\begin{pspicture}(-5,-3)(5,1)
\pnode(3;60){B}\pnode(3;300){A}\pnode(6;0){V}\pnode(4;0){K}
\pspolygon(0,0)(A)(V)(B) \psline[linewidth=0.5pt](2.5;180)(V)
\pscurve(B)(K)(A) \psline[linestyle=dashed](K)(B)
\uput[dr](K){$K$}\uput[dl](0,0){$P$}\uput[60](B){$B$}
\uput[270](1.7;0){$x_K$}\uput{0.5}[90](1.2;60){$r_0$}
\end{pspicture}
 
\caption{Existence under local symmetry (Prop. 4.2)}\end{center}
\end{figure}

\emph{Proof.} (Fig. 14.) $u(\pi/6)>0$ and $u(7\pi/6)>0$ follow
from symmetry exactly as before, so we only need to check
$w(\pi/6)<0$. Letting $P,\bar P$ be the points on the chord
$\bar{K}K$ where the normals through $B(\pi/6), B(5\pi/6)$ (resp.
through $B(7\pi/6),B(11\pi/6)$) meet, this is equivalent to
$|KP|+|\bar{K}\bar{P}|<L_{chord}$; so let $x=x_K=|KP|$; we
estimate $x$ in terms of $r_K$. First note (Fig.14):
$$\frac{\pi}3=\int_K^B\frac 1r ds\geq \frac 1{r_K}arclength_{\cal
C}(KB)\geq \frac 1{r_K}|KB|.$$ Then with $r_0:=|BP|=-\langle
B-P,N(5\pi/6)\rangle$ and
$B=P+(x,0)+\int_{\pi/2}^{5\pi/6}r(\tau)T(\tau)d\tau$, we have:
$$r_0=\frac x2+\int_{\pi/2}^{5\pi/6}r(\tau)\sin(\frac{5\pi}6-\tau)d\tau,$$
so $\frac x2<r_0<\frac 12 (x+r_K)$ and (using the triangle $BPK$,
where the interior angle at $P$ is $\pi/3$):
$$r_K^2\frac{\pi^2}9\geq |KB|^2=r_0^2+x^2-xr_0\geq
\frac{3x^2}4-\frac{r_K}2x,$$ which gives the upper bound: $x\leq
c_0 r_K$. Repeating the argument at the other end for
$x_{\bar{K}}=|\bar{K}\bar{P}|$, we get $x_{\bar{K}}\leq
c_0r_{\bar{K}}$, so:
$$x_K+x_{\bar K}\leq c_0(r_K+r_{\bar K})<L_{chord}.$$
\vspace{.2cm}

Our last existence result in this section makes no symmetry
assumptions, but yields a weaker conclusion.\vspace{.2cm}

\textbf{Proposition 4.3.} \emph{Let $\cal C$ be a strictly convex
curve (defined by $p\in {\cal P}_{pw}$). Assume $\cal C$ has only
\emph{two} critical chords (both non-degenerate), making an angle
greater than $\pi/3$ (unoriented angles, taking values in $(0,
\pi/2]$). Then the parallel curves ${\cal C}_{\epsilon}$ support a
double triode stably, for all $\epsilon >0$ sufficiently large
(and possibly for $\epsilon=0$ already).}\vspace{.2cm}

\emph{Proof.} We only need to find $\theta \in [0,\pi)$ so that
$y(\theta)=0$ and $w(\theta)<0$; then add a sufficiently large
constant to $p$ to guarantee $u(\theta),u(\bar{\theta})$ are
positive (adding constants to $p$ leaves $w$ and $y$ unchanged.)
Let the minimal chord be $B(0)B(\pi)$ and the maximal
$B(\alpha)B(\alpha+\pi)$, with $\alpha \in (\pi/3, \pi/2]$. Then
$w<0$ in $(0,\alpha)\cup (\pi, \alpha +\pi)$ and $w>0$ in
$(\alpha,\pi)\cup (\alpha+\pi,2\pi)$. Letting $\tilde
w(\theta):=w(\theta+2\pi/3)$, this means: $\tilde{w}<0$ in
$(\alpha+\pi,\alpha+\frac{4\pi}3)\cup
(\frac{\pi}3,\alpha+\frac{4\pi}3)$, $\tilde{w}>0$ in
$(\alpha+\frac{\pi}3, \frac{4\pi}3)\cup
(\alpha+\frac{4\pi}3,\frac{7\pi}3).$ Since $\alpha \in
(\frac{\pi}3,\frac{\pi}2]$, there are arcs in $S$ where both $w$
and $\tilde{w}$ are positive, as well as arcs where both are
negative, determining the sign of $y$:
$$y=w+\tilde{w}<0\mbox{ in }(\frac{\pi}3,\alpha)\cup
(\frac{4\pi}3, \alpha +\pi);$$
$$y=w+\tilde{w}>0\mbox{ in }(\alpha+\frac{\pi}3,\pi)\cup
(\alpha+\frac{4\pi}3,2\pi).$$ Thus $y$ must have a zero $\theta
\in (0,\frac{\pi}3)$, where $w<0$. (There is also a zero in the
interval $(\alpha, \alpha+\frac{\pi}3)$, where $w>0$, but we'll
ignore it.) \vspace{.2cm}

\emph{Remark 4.1.} We can guarantee $u>0$ already for $\cal C$ by
imposing a `curvature pinching condition', as in Proposition 3.2
(e.g. $(1+\sqrt{3}/2)\min_{\theta}r>\max_{\theta}r$.)\vspace{.2cm}

\emph{Remark 4.2.} The restriction to two critical chords is made
just to simplify the statement; a more general result follows from
the same argument (assuming all critical chords are non-degenerate),
but is cumbersome to state. Consider the \emph{oriented} angular
distance from a \emph{minimal} chord to the next \emph{maximal}
chord (moving counterclockwise on $S$), with values in $(0,\pi)$. If
this distance is always greater than $2\pi/3$, the same conclusion
follows. Another sufficient condition is: the smallest oriented
distance from a \emph{maximal} chord to the next \emph{minimal}
chord is less than $2\pi/3$, \emph{and} if chords achieving this
least distance define the arc $(\theta_{max},\theta_{min})$ in $S$
(where $w>0$), the union of this arc with the two adjacent arcs
where $w<0$ has angular measure at least $2\pi/3$. We omit the
proof. \vspace{.3cm}

\emph{4.3 Examples without double triodes.}\vspace{.2cm}

The critical points of parallelogram length are solutions of
$y(\theta):=w(\theta)+w(\theta+2\pi/3)=0$. To define a double
triode, we must in addition have $w(\theta)<0$. Thus if at every
zero of $y$ (in a fundamental domain, such as $[0,\pi)$) we have
$w>0$ (or, equivalently, $y(\theta+\pi/3)>y(\theta+2\pi/3)$), then
$\cal C$ does not support a double triode (nor does any of its
outer parallel curves.)

The idea to construct examples is to start from $p_0$ of constant
width ($w_0\equiv 0$; e.g., the circle), then perturb it by adding
$\tilde{p}$ (not necessarily convex!), constructed so that
$\tilde{y}$ has the desired property. This leads to the following
`instability property' for curves of constant width. \vspace{.2cm}

\textbf{Proposition 4.4.}\emph{ Let ${\cal C}_0$ (with support
function $p_0\in {\cal P}_{smooth}$) be a strictly convex curve of
constant width. Then arbitrarily $C^2$ close to $p_0$ one finds
(i) support functions $p$ so that (stably) neither the convex
curve $\cal C$ nor its outer parallel curves support double
triodes; (ii) support functions $p$ of curves for which (stably)
some outer parallel curve carries double triodes.} \vspace{.2cm}

\emph{Proof.} Let $\tilde{w}$ ($\pi$-periodic, with
$\int_0^{\pi}\tilde{w}(\tau)d\tau=0$) have the property:
$\tilde{w}>0$ at each zero of $\tilde{w}(.)+\tilde{w}(.+2\pi/3)$,
and these zeros are all transversal. Find $\tilde{p}$
$2\pi$-periodic so that
$\tilde{w}(\theta)=-(\tilde{p}(\theta)+\tilde{p}(\theta+\pi))$
(for example, $\tilde{p}=(-1/2)\int^{\theta}\tilde{w}$). Then the
support function:
$$p(\theta)=\frac
1{\lambda+1}\tilde{p}(\theta)+\frac{\lambda}{\lambda+1}p_0(\theta)$$
has, for all $\lambda>0$ large enough, the properties: (a) $p$ is
positive and strictly convex (meaning $r>0$); (b) neither $p$ nor
$p+c$ (for any $c>0$) supports double triodes, since $w>0$ at each
zero of $w(.)+w(.+2\pi/3)$ (given that this is true for
$\tilde{w}$ and $w=\frac 1{\lambda+1}\tilde{w}$). Since $p$ can be
made arbitrarily close to $p_0$ in $C^2$ norm, this proves (i). To
show (ii), consider the support function:
$$\bar{p}(\theta)=-\frac
1{\lambda+1}\tilde{p}(\theta)+\frac{\lambda}{\lambda+1}p_0(\theta),$$
for which $\bar{w}=-(\lambda+1)^{-1}\tilde{w}$, so $\bar{w}<0$ at
each zero of $\bar{w}(.)+\bar{w}(.+2\pi/3)$. Again, for large
enough $\lambda$, $\bar{p}$ is positive and $\bar{r}>0$, so
$\bar{p}$ defines a strictly convex curves which supports double
triodes, possibly after adding a sufficiently large constant (this
is not needed if $u_0>0$ everywhere.)\vspace{.2cm}

It remains to exhibit a $\pi$-periodic function $\tilde{w}$ with
the desired properties. The $2\pi$-periodic function:
$$z(\theta)=\frac{M-M\cos(\theta-\bar{\theta})-(1-\cos
\theta)(1-\cos \bar{\theta})}{M-M\cos(\theta-\bar{\theta})+(1-\cos
\theta)(1-\cos \bar{\theta})}$$ (where $M$ and $\bar{\theta}$ are
parameters) has $0$ and $\bar{\theta}$ as its only critical points
($z(0)=1, z(\bar{\theta})=-1$ and $-1<z(\theta)<1$ otherwise).
Choosing $M=1,\bar{\theta}=\pi/6$, one finds that both at
$\theta=0$ and $\theta=\bar{\theta}$:
$$z'(\theta+\frac{2\pi}3)>z'(\theta+\frac{4\pi}3).$$
Let $Y(\theta)=z(2\theta)$. Then $0$ and $\bar{\theta}/2=\pi/12$
are the only critical points of $Y$ in $[0,\pi)$, and at each of
them: $Y'(\theta+\frac{\pi}3)>Y'(\theta+\frac{2\pi}3)$; thus
setting $\tilde{y}=Y'$ and
$\tilde{w}(\theta)=(1/2)(\tilde{y}(\theta)+\tilde{y}(\theta+\pi/3)-\tilde{y}(\theta+2\pi/3))$,
we have $\tilde{w}>0$ when $\tilde{y}=0$, as desired.
\vspace{.2cm}

\emph{Example 4.1.} The construction can be made completely
explicit. For the given values $M=1$, $\bar{\theta}=\pi/6$, one
finds $\tilde{w}(0)\sim 0.0635$ and $\tilde{w}(\pi/12)\sim 0.0125$
(both positive). We may take
$\tilde{p}(\theta)=(-1/2)(Y(\theta)+Y(\theta+\pi/3)-Y(\theta+2\pi/3))$
(then $\tilde{w}=-\tilde{p}'$).

Adding a large constant to this $\tilde{p}$ (which corresponds to
perturbing a circle by $\tilde{p}$), we get an explicit example of
a strictly convex curve supporting no double triodes (one finds
numerically that $\tilde{p}+300$ already corresponds to a convex
curve). Or one can start from the smoothed Reuleaux triangle
($p_{\epsilon}$ in (1.4)) and let:
$$p=\frac
1{\lambda+1}\tilde{p}+\frac{\lambda}{\lambda+1}p_{\epsilon},\quad
\bar{p}=-\frac
1{\lambda+1}\tilde{p}+\frac{\lambda}{\lambda+1}p_{\epsilon}.$$ For
$\epsilon=1$ and $\lambda=1,000$, both represent convex curves;
$\bar{p}$ stably supports double triodes, while $p$ stably
doesn't- and both can be taken arbitrarily close to
$p_{\epsilon}$. \vspace{.2cm}

\emph{Remark 4.3.} For future reference, we record here the
following non-existence criterion: if $p\in {\cal P}_{pw}$ has the
property that, for all $\theta$ so that $y(\theta)=0$, either
$w(\theta)>0$ or $u(\theta)<0$, then $\cal C$ supports no double
triodes (stably.)\vspace{.2cm}

\emph{Remark 4.4.} In example 4.1, we could take $p_{\epsilon}$ to
have constant width and height (say, a finite linear combination
of Fourier components with frequency an odd integer not divisible
by 3.) With a little more care, the perturbation $\tilde{p}$ can
also be chosen to have vanishing height function; this gives a
construction of examples of convex constant-height curves which do
not support double triodes (however, all known examples do support
triodes.)


\newpage

\textbf{5. Hexagonal cells.}\vspace{.2cm}

\emph{5.1 The hexagonal cell configuration; holonomy.} An
hexagonal cell in $D$ (when it exists) is defined by a choice of
six points on $\cal C$ (with angular separation $\pi/3$), and by
an additional parameter taking values in an open interval,
measuring how far the cell is from the boundary. The existence of
a cell for a given $\theta\in S$ is determined by: (i) a
criticality condition, depending only on the configuration of
normal lines defined by $\theta$ and its $\pi/3$ translates (i.e.,
unchanged for parallel curves); (ii) inequalities ruling out
certain critical configurations of normal lines (analogous to
$w<0$ for double triodes); (iii) requiring points of the
configuration to be inside $\cal C$, which can always be achieved
by taking outer parallel curves. \vspace{.2cm}

The configuration of tangent polygons and normal lines to be
considered for hexagons \emph{includes} those considered for
triodes and double triodes, and it is useful to preserve the
notation used in those two cases. We describe the \emph{notational
conventions} used in this section (see Fig.15.)

(1) Indices 1,2,3 denote positive $2\pi/3$ shifts from $\theta$,
and a bar is used for $\pi$ shifts. Thus the boundary points are:
$B_1, \bar{B_3},B_2,\bar{B_1},B_3,\bar{B_2}$ (in cyclic order,
with $B_1=B(\theta)$), and the unit tangent and normal vector at
these points are denoted accordingly.

(2) In general, the \emph{cyclic order} of a generic `index' $I$
taking 6 values will be: $I=(1,\bar{3},2,\bar{1},3,\bar{2})$.
Advancing one step in this cyclic order corresponds to shifting
$\theta$ by $\pi/3$.

(3) From the circumscribed equilateral triangles ${\cal
T}(\theta),{\cal T}(\bar{\theta})$ we have the 6 `triangular
normal intersections' $P_1,P_2,P_3,\bar{P_1},\bar{P_2},\bar{P_3}$.
(These are indicated by dots in Fig. 15; to avoid encumbering the
figure, only $P_1$ and $\bar{P_2}$ are labelled.) The derived
height $\omega$ and `triangle functions' $u,v$ are defined exactly
as before (and e.g.
$\bar{u}_2=u(\theta+\frac{2\pi}3+\pi)=u(\theta+\frac{5\pi}3)$).

(4) In addition to the `projections' of the $P_i$ considered in
section 4: $$U_1=pr_{\bar{2}1}(\bar{P_1}),\quad
\bar{U_1}=pr_{2\bar{1}}(P_1)$$ (and cyclic, e.g.
$\bar{U_3}=pr_{1\bar{3}}(P_3)$), we also need `backward'
projections such as:
$$V_1=pr_{\bar{n}_3\rightarrow
n_1}(\bar{P_3})=pr_{\bar{3}1}(\bar{P_3})$$ (and cyclic:
$\bar{V_3}=pr_{2\bar{3}}(P_2),V_2=pr_{\bar{1}2}(\bar{P_1})$, etc.
).

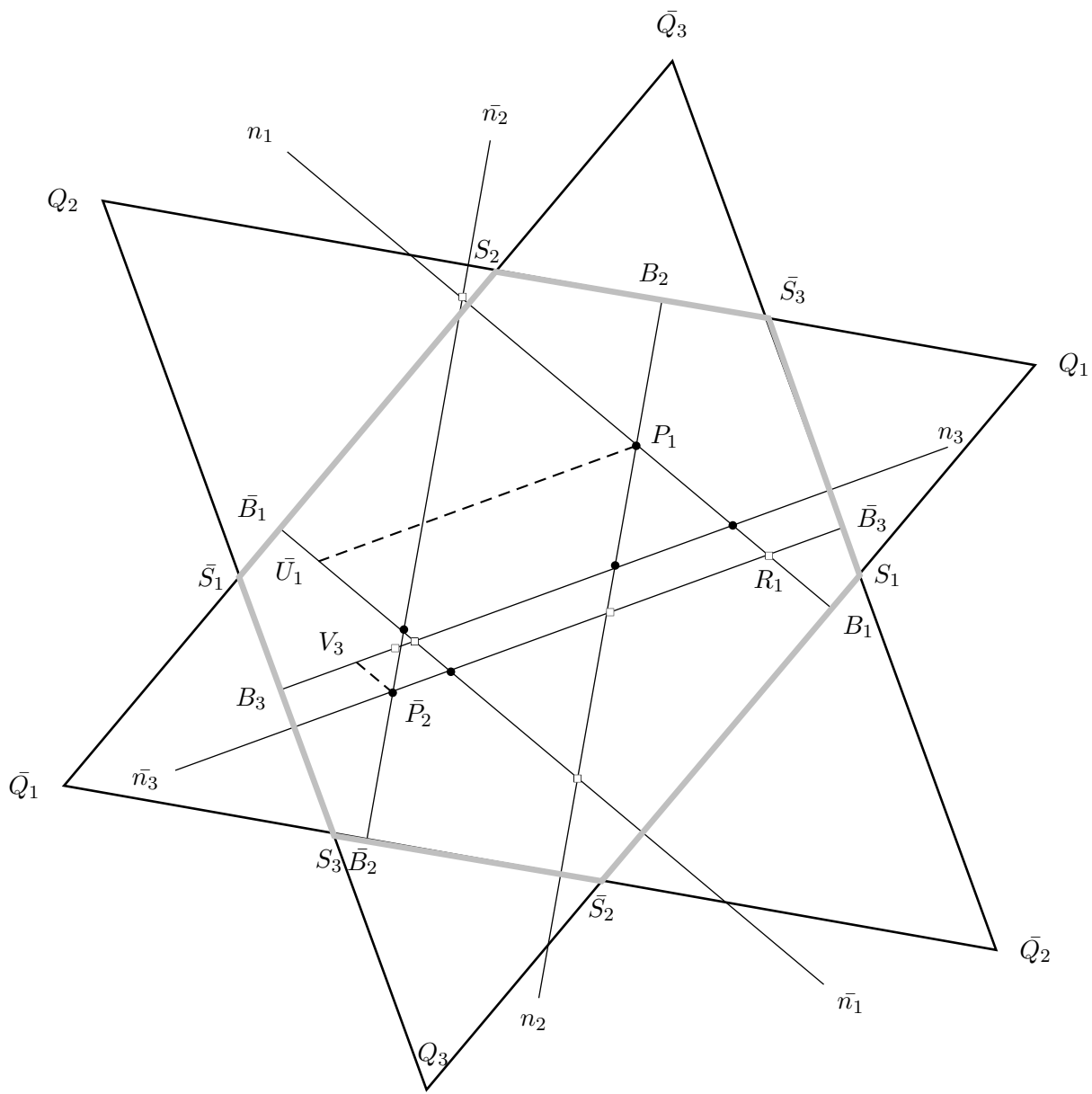
\begin{figure} \psset{unit=0.7cm}
\begin{center}
\SpecialCoor
\begin{pspicture}(-3,0)(7,16)
\pnode(0,0){Q3}  \pnode([angle=50,nodesep=20]Q3){Q1}
\pnode([angle=110,nodesep=20]Q3){Q2}
\pspolygon[linewidth=1pt](Q3)(Q1)(Q2)
\pnode([angle=50,nodesep=13.3]Q3){B1}
\pnode([angle=110,nodesep=9]Q3){B3}
\pnode([angle=170,nodesep=8]Q1){B2}
\psline[linewidth=0.5pt](B1)([angle=140,nodesep=15]B1)
\uput[140]([angle=140,nodesep=15.1]B1){$n_1$}
\psline[linewidth=0.5pt](B2)([angle=260,nodesep=15]B2)
\uput[260]([angle=260,nodesep=15.1]B2){$n_2$}
\psline[linewidth=0.5pt](B3)([angle=20,nodesep=15]B3)
\uput[20]([angle=21,nodesep=14.6]B3){$n_3$}
\pnode(10;140){Q1bar} \pnode([angle=350,nodesep=20]Q1bar){Q2bar}
\pnode([angle=50,nodesep=20]Q1bar){Q3bar}
\pspolygon[linewidth=1pt](Q1bar)(Q2bar)(Q3bar)
\pnode([angle=350,nodesep=6.5]Q1bar){B2bar}
\pnode([angle=50,nodesep=7.1]Q1bar){B1bar}
\pnode([angle=110,nodesep=9.5]Q2bar){B3bar}
\psline[linewidth=0.5pt]([angle=80,nodesep=15]B2bar)(B2bar)
\uput[80]([angle=80,nodesep=15.1]B2bar){$\bar{n_2}$}
\psline[linewidth=0.5pt]([angle=320,nodesep=15]B1bar)(B1bar)
\uput[320]([angle=320,nodesep=15.1]B1bar){$\bar{n_1}$}
\psline[linewidth=0.5pt]([angle=200,nodesep=15]B3bar)(B3bar)
\uput[200]([angle=200,nodesep=15.1]B3bar){$\bar{n_3}$}
\pnode([angle=80,nodesep=3.143]B2bar){P2bar}
\pnode([angle=260,nodesep=3.143]B2){P1}
\pnode([angle=260,nodesep=5.714]B2){P2}
\pnode([angle=140,nodesep=2.714]B1){P3}
\pnode([angle=80,nodesep=4.5]B2bar){P1bar}
\pnode([angle=320,nodesep=4.714]B1bar){P3bar}
\psdots(P1)(P2)(P3)(P1bar)(P2bar)(P3bar)
\pnode([angle=140,nodesep=1.714]B1){R1}
\pnode([angle=200,nodesep=5.2214]B3bar){R3bar}
\pnode([angle=260,nodesep=10.286]B2){R2}
\pnode([angle=320,nodesep=3.714]B1bar){R1bar}
\pnode([angle=20,nodesep=2.5714]B3){R3}
\pnode([angle=80,nodesep=11.643]B2bar){R2bar}
\psdots[dotstyle=square](R1)(R2)(R3)(R1bar)(R2bar)(R3bar)
\pnode([angle=140,nodesep=1]P2bar){V3}
\pnode([angle=200,nodesep=7.143]P1){U1bar}
\psline[linestyle=dashed](P1)(U1bar)
\psline[linestyle=dashed](P2bar)(V3)
\pnode([angle=50,nodesep=5.757]Q3){S2bar}
\pnode([angle=110,nodesep=5.714]Q3){S3}
\pnode([angle=170,nodesep=5.714]Q1){S3bar}
\pnode([angle=50,nodesep=5.757]Q1bar){S1bar}
\pnode([angle=50,nodesep=14.186]Q1bar){S2}
\pnode([angle=110,nodesep=8.429]Q2bar){S1}
\pspolygon[linewidth=2.5pt,linecolor=lightgray](S2bar)(S1)(S3bar)(S2)(S1bar)(S3)
\uput{0.5}[80](Q3){$Q_3$} \uput{0.5}[0](Q1){$Q_1$}
\uput{0.5}[180](Q2){$Q_2$}\uput{0.3}[320](B1){$B_1$}
\uput{0.3}[105](B2){$B_2$}
\uput{0.3}[200](B3){$B_3$} \uput{0.5}[180](Q1bar){$\bar{Q_1}$}
\uput{0.5}[0](Q2bar){$\bar{Q_2}$}
\uput{0.5}[90](Q3bar){$\bar{Q_3}$}
\uput{0.3}[140](B1bar){$\bar{B_1}$}
\uput{0.2}[260](B2bar){$\bar{B_2}$}
\uput{0.3}[20](B3bar){$\bar{B_3}$} \uput{0.3}[20](P1){$P_1$}
\uput{0.3}[320](P2bar){$\bar{P_2}$} \uput{0.3}[270](R1){$R_1$}
\uput{0.3}[200](U1bar){$\bar{U_1}$} \uput{0.3}[140](V3){$V_3$}
\uput{0.3}[260](S3){$S_3$} \uput{0.3}[270](S2bar){$\bar{S_2}$}
\uput{0.3}[0](S1){$S_1$} \uput{0.4}[50](S3bar){$\bar{S_3}$}
\uput{0.2}[110](S2){$S_2$} \uput{0.3}[180](S1bar){$\bar{S_1}$}


\end{pspicture}
 
\caption{The configuration for an hexagonal cell}\end{center}
\end{figure}

Recall $U_I=B_I+(2/\sqrt{3})\sigma_IN_I$,where in (4.1) we found
an expression for $\sigma_I$. In the same way we let:
$$V_1=B_1+\frac 2{\sqrt{3}}\tau_1N_1,$$
with $\tau_1=\tau(\theta)$, and use the fact that also
$V_1=\bar{P}_3+\beta N_2$ (for some $\beta \in \mathbb{R}$) to
compute $\tau_1$, by taking inner products with $T_2$. We find:
$$\tau_1=\langle \bar{B}_3-B_1,T_2\rangle +\bar{u}_3.$$
Using $\bar{u}_3=\langle \bar{B}_1-\bar{B}_3,\bar{T}_1\rangle$ and
rearranging:
$$\tau_1=\langle B_1-B_3,T_3\rangle-\langle
\bar{B}_3-B_3,T_3\rangle-\langle \bar{B}_1-B_1,T_1\rangle,$$ which
can also be written in the form:
$$\tau(\theta)=v(\theta)-w(\theta)-w(\theta+\frac{4\pi}3)\quad
(5.1)$$

Note that there is a total of 12 points defined by both types of
projection (6 of type $U_I$, 6 of type $V_I$). (Only $\bar{U_1}$
and $V_3$ are shown in Fig.15.)

(5) The configuration includes also 6 new points, `hexagonal
normal intersections': $R_1,\bar{R}_3,R_2,\bar{R}_1,R_3,\bar{R}_2$
(in cyclic order), where: $\{R_1\}=n_1\cap \bar{n}_3$ and cyclic
(e.g. $\{\bar{R}_3\}=\bar{n}_3\cap n_2$,
$\{\bar{R}_2\}=\bar{n}_2\cap n_1$, etc.) These points are
indicated by white squares in Fig.15, where only $R_1$ is
labelled. Define functions $s(\theta)=s_1, t(\theta)=t_1$ by:
$$R_1=B_1+\frac2{\sqrt{3}}s_1N_1,\quad \bar{R}_2=B_1+\frac
2{\sqrt{3}}t_1N_1.$$

We compute $s_1$ and $t_1$ in the usual way and find:
$$s_1=\langle \bar{B}_3-B_1,\bar{T}_3\rangle, \quad
s(\theta)=\langle B(\theta
+\frac{\pi}3)-B(\theta),T(\theta+\frac{\pi}3)\rangle,$$
$$t_1=\langle B_1-\bar{B}_2,\bar{T}_2\rangle,\quad
t(\theta)=\langle
B(\theta)-B(\theta-\frac{\pi}3),T(\theta-\frac{\pi}3)\rangle.$$ It
turns out $s$ and $t$ can be expressed in terms of previously
defined functions. We have:
$$s_1=\langle \bar{B}_3-B_3,\bar{T}_3\rangle+\langle
B_3-B_1,\bar{T}_3\rangle=-w_3+\langle
B_1-B_3,T_3\rangle=-w_3+v_1,$$ or:
$$s(\theta)=v(\theta)-w(\theta+\frac{4\pi}3)\qquad (5.2)$$
A similar computation yields:
$$t(\theta)=u(\theta)+w(\theta+\frac{2\pi}3)\qquad (5.3)$$
\vspace{.2cm}

In total, the configuration for hexagonal cells includes 24
intersection points, defined by the configuration of normals (at
$\theta$ and its $\pi/3$ translates.) \vspace{.2cm}

\emph{Remark 5.1.} From (5.3) and expression (4.1) for $\sigma$,
we find:
$$\sigma(\theta)=t(\theta)+w(\theta)\qquad (5.4),$$
and combining (5.2) and (5.1), we similarly find:
$$\tau(\theta)=s(\theta)-w(\theta)\qquad (5.5).$$
In particular, the functions $s,t,\sigma$ and $\tau$ can be written
as simple combinations of $u$, $\omega$ and their translates:
$$\begin{array}{rcl}
s&=&-\omega(\theta)-w(\theta+\frac{\pi}3)+u(\theta);\\
t&=&w(\theta+\frac{2\pi}3)+u(\theta);\\
\sigma&=&w(\theta)+w(\theta+\frac{2\pi}3)+u(\theta);\\
\tau&=&-\omega(\theta)-w(\theta)-w(\theta+\frac{\pi}3)+u(\theta).
\end{array}\qquad (5.6)$$
We see immediately that by adding a positive constant to the
support function $p(\theta)$ (which adds a constant to $u$ without
changing the other functions)-i.e., by considering outer parallel
curves- we can make all these functions positive.\vspace{.2cm}

The intersections of the equilateral triangles ${\cal
T}(\theta)=(Q_1,Q_2,Q_3),{\cal
T}(\bar{\theta})=(\bar{Q}_1,\bar{Q}_2,\bar{Q}_3)$ (for a given
$\theta$) define six points
($S_1,\bar{S}_3,S_2,\bar{S}_1,S_3,\bar{S}_2$), vertices of a
circumscribed equiangular hexagon ${\cal H}(\theta)$. To construct
a hexagonal cell ${\mathbb
H}(\theta)=(H_1\bar{H}_3H_2\bar{H}_1H_3\bar{H}_2)$, we pick a
small positive $x=x_1$ and define:
$$H_1=B_1+\frac{2}{\sqrt{3}}x_1N_1,$$
and successively:
$$\bar{H}_3=pr_{1\bar{3}}(H_1)=\bar{B}_3+\frac
2{\sqrt{3}}\bar{x}_3 \bar{N}_3,\ldots, H_{I+1}=pr_{n_I\rightarrow
n_{I+1}}(H_I),\ldots$$ and finally:
$$H_1^*=pr_{\bar{n}_2\rightarrow n_1}(\bar{H}_2)=B_1+\frac
2{\sqrt{3}}x_1^*N_1.$$ In general $H_1^*\neq H_1$, and instead of
a `cell' we have only a hexagonal `chain' with seven vertices,
beginning with $H_1$ and ending with $H_1^*$, both on $n_1$. (Fig.
16; points $\bar{H_3}$, $H_3$ and $H_1^*$ are not labelled, to
avoid overloading the figure.) We call this `defect' the
\emph{holonomy} of $\theta$:
$$hol(\theta)=x_1^*-x_1;$$
we will soon see that it depends only on the configuration of
normals. Its vanishing is the `criticality condition' mentioned
above. To compute the holonomy, beginning with:
$$\langle R_1-H_1,N_1\rangle=\langle
R_1-\bar{H}_3,\bar{N}_3\rangle$$ we obtain:
$$\begin{array}{rcl}
\bar{x}_3&=& x_1+\bar{t}_3-s_1,\mbox{ and successively in cyclic
order: }\\
x_2&=&\bar{x}_3+t_2-\bar{s}_3,\ldots\\
x_1^*&=&\bar{x}_2+t_1-\bar{s}_2. \end{array}$$ Adding the results,
we have:
$$hol(\theta)=x_1^*-x_1=\sum_{i=1}^3(t_i-s_i)+\sum_{i=1}^3
(\bar{t}_i-\bar{s}_i).$$ Using the expressions given above for
$s_I$ and $t_I$, one easily computes that:
$$\begin{array}{rcl}
hol(\theta)&=&\langle \bar{B}_1-B_1,T_1\rangle +\langle
\bar{B}_2-B_2,T_2\rangle+\langle \bar{B}_3-B_3,T_3\rangle\\
&=&w(\theta)+w(\theta+2\pi/3)+w(\theta+4\pi/3).\end{array}$$

\begin{figure} \psset{unit=0.5cm}
\begin{center}
\SpecialCoor
\begin{pspicture}(0,3)(8,18)
\pnode(0,0){Q3}  \pnode([angle=50,nodesep=20]Q3){Q1}
\pnode([angle=110,nodesep=20]Q3){Q2}
\pnode([angle=50,nodesep=13.3]Q3){B1}
\pnode([angle=110,nodesep=9]Q3){B3}
\pnode([angle=170,nodesep=8]Q1){B2}
\psline[linewidth=0.5pt](B1)([angle=140,nodesep=11]B1)
\uput[140]([angle=140,nodesep=11.1]B1){$n_1$}
\psline[linewidth=0.5pt](B2)([angle=260,nodesep=13]B2)
\uput[260]([angle=260,nodesep=13.1]B2){$n_2$}
\psline[linewidth=0.5pt](B3)([angle=20,nodesep=15]B3)
\uput[20]([angle=21,nodesep=14.6]B3){$n_3$}
\pnode(10;140){Q1bar} \pnode([angle=350,nodesep=20]Q1bar){Q2bar}
\pnode([angle=50,nodesep=20]Q1bar){Q3bar}
\pnode([angle=350,nodesep=6.5]Q1bar){B2bar}
\pnode([angle=50,nodesep=7.1]Q1bar){B1bar}
\pnode([angle=110,nodesep=9.5]Q2bar){B3bar}
\psline[linewidth=0.5pt]([angle=80,nodesep=13]B2bar)(B2bar)
\uput[80]([angle=80,nodesep=13.1]B2bar){$\bar{n_2}$}
\psline[linewidth=0.5pt]([angle=320,nodesep=11]B1bar)(B1bar)
\uput[320]([angle=320,nodesep=11.1]B1bar){$\bar{n_1}$}
\psline[linewidth=0.5pt]([angle=200,nodesep=15]B3bar)(B3bar)
\uput[200]([angle=200,nodesep=15.1]B3bar){$\bar{n_3}$}
\pnode([angle=80,nodesep=3.143]B2bar){P2bar}
\pnode([angle=260,nodesep=3.143]B2){P1}
\pnode([angle=260,nodesep=5.714]B2){P2}
\pnode([angle=140,nodesep=2.714]B1){P3}
\pnode([angle=80,nodesep=4.5]B2bar){P1bar}
\pnode([angle=320,nodesep=4.714]B1bar){P3bar}
\pnode([angle=140,nodesep=1.714]B1){R1}
\pnode([angle=200,nodesep=5.2214]B3bar){R3bar}
\pnode([angle=260,nodesep=10.286]B2){R2}
\pnode([angle=320,nodesep=3.714]B1bar){R1bar}
\pnode([angle=20,nodesep=2.5714]B3){R3}
\pnode([angle=80,nodesep=11.643]B2bar){R2bar}
\pnode([angle=140,nodesep=1]P2bar){V3}
\pnode([angle=200,nodesep=7.143]P1){U1bar}
\pnode([angle=260,nodesep=2.6]U1bar){W1}
\pnode([angle=80,nodesep=1.3]V3){W3}
\pnode([angle=50,nodesep=5.757]Q3){S2bar}
\pnode([angle=110,nodesep=5.714]Q3){S3}
\pnode([angle=170,nodesep=5.714]Q1){S3bar}
\pnode([angle=50,nodesep=5.757]Q1bar){S1bar}
\pnode([angle=50,nodesep=14.186]Q1bar){S2}
\pnode([angle=110,nodesep=8.429]Q2bar){S1}
\psccurve[linewidth=1pt](B3)(B2bar)(B1)(B3bar)(B2)(B1bar)
\pnode([angle=140,nodesep=1.5]B1){H1} \uput{0.1}[260](H1){$H_1$}
\pnode([angle=80,nodesep=0.2]H1){H3bar}
\pnode([angle=140,nodesep=3.9]H3bar){H2}
\uput{0.25}[20](H2){$H_2$}
\pnode([angle=200,nodesep=7.3]H2){H1bar}\uput{0.1}[80](H1bar){$\bar{H_1}$}
\pnode([angle=260,nodesep=2.9]H1bar){H3}
\pnode([angle=320,nodesep=2.6]H3){H2bar}\uput{0.1}[200](H2bar){$\bar{H_2}$}
\pnode([angle=20,nodesep=10]H2bar){H1star} \psline[linecolor=blue]
{->}(H1)(H3bar)(H2)(H1bar)(H3)(H2bar)(H1star)
\uput{0.3}[320](B1){$B_1$}
\uput{0.3}[105](B2){$B_2$}
\uput{0.3}[200](B3){$B_3$} 
\uput{0.3}[140](B1bar){$\bar{B_1}$}
\uput{0.2}[260](B2bar){$\bar{B_2}$}
\uput{0.3}[20](B3bar){$\bar{B_3}$} 

\end{pspicture}
 
\caption{An hexagonal chain with holonomy}
\end{center}
\end{figure}
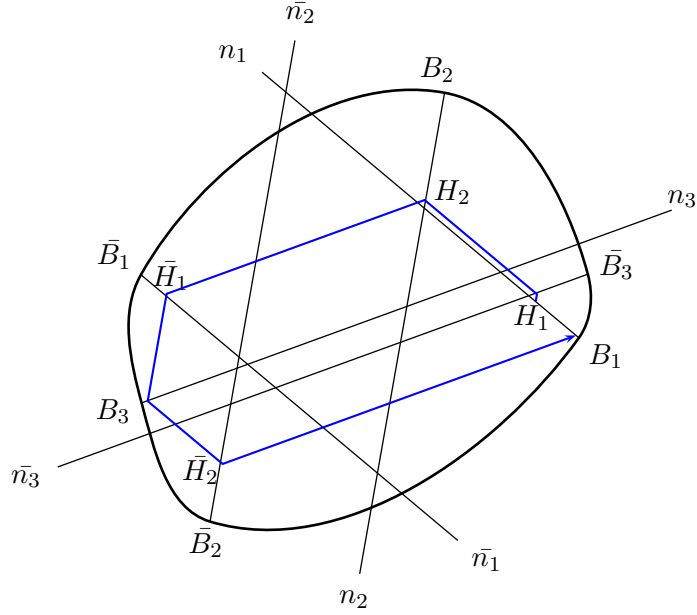

\emph{Remark 5.2.} Recalling that
$w(\theta)=-(p'(\theta)+p'(\bar{\theta}))$, we have the equivalent
expression:
$$hol(\theta)=-\sum_{j=0}^5p'(\theta+j\pi/3)
=-(\omega(\theta)+\omega(\bar{\theta})).$$ \vspace{.2cm}

We now relate the holonomy to the perimeter $L_{hex}(\theta)$ of
the circumscribed hexagon ${\cal H}(\theta)$. Proceeding as usual,
we find:
$$\begin{array}{rcl}
S_1&=&B_1+\frac 2{\sqrt{3}}\langle B_1-\bar{B}_3,\bar{N}_3\rangle
T_1=\bar{B}_3+\frac 2{\sqrt{3}}\langle B_1-\bar{B}_3,N_1\rangle
\bar{T}_3,\\
\bar{S}_3&=&\bar{B}_3+\frac 2{\sqrt{3}}\langle
\bar{B}_3-B_2,N_2\rangle \bar{T}_3, \end{array}$$ which gives for
the length of the side $S_1\bar{S}_3$:
$$\langle \bar{S_3}-S_1,\bar{T}_3\rangle=-\frac
2{\sqrt{3}}(\langle B_1,N_1\rangle+\langle B_2,N_2\rangle+\langle
B_3,N_3\rangle).$$ Repeating this for the other sides and adding
the results (using $\bar{N}_i=-N_i$), we obtain:
$$\begin{array}{rcl}
L_{hex}(\theta)&=&\frac 2{\sqrt{3}}(\langle
\bar{B}_1-B_1,N_1\rangle+\langle \bar{B}_2-B_2, N_2\rangle+\langle
\bar{B}_3-B_3, N_3\rangle)\\
&=&\frac
2{\sqrt{3}}(b(\theta)+b(\theta+\frac{2\pi}3)+b(\theta+\frac{4\pi}3)).
\end{array}$$
This immediately implies the following (since $w=-b'$):
\vspace{.2cm}

\textbf{Proposition 5.1}
\emph{$hol(\theta)=-(\sqrt{3}/2)L_{hex}'(\theta)$. Hence there are
always at least two directions $\theta\in S$ along which the
`holonomy' vanishes.} \vspace{.2cm}

\emph{Remark 5.3.} We can also relate the hexagonal perimeter and
the height function: from
$h(\theta)=p(\theta)+p(\theta+2\pi/3)+p(\theta+4\pi/3)$ and
$b(\theta)=p(\theta)+p(\bar{\theta})$ follows:
$$L_{hex}(\theta)=\frac 2{\sqrt{3}}(h(\theta)+h(\bar{\theta})).$$
\vspace{.2cm}

\emph{Remark 5.4 (Total length of the network.)} The total length
of a closed cell ${\mathbb
H}=(H_1\bar{H}_3H_2\bar{H}_1H_3\bar{H}_2)$ is:
$$\begin{array}{rcl}L(\theta)&=&\langle
\bar{H}_3-H_1,\bar{N}_2\rangle+\langle H_2-\bar{H}_3,N_1\rangle
+\langle \bar{H}_1-H_2,\bar{N}_3\rangle\\&+&\langle
H_3-\bar{H}_1,N_2\rangle +\langle
\bar{H}_2-H_3,\bar{N}_1\rangle+\langle H_1-\bar{H}_2,N_2\rangle\\
&+&\frac 2{\sqrt{3}}(x_1+\bar{x}_3+\ldots+\bar{x}_2).\end{array}$$
We rewrite the terms on the right-hand side in the form:
$$
\langle \bar{H}_3-H_1,\bar{N}_2\rangle= \langle
\bar{B_3}-B_1,\bar{N}_2\rangle -\frac 1{\sqrt{3}}(x_1+\bar{x}_3)$$
(and cyclic), and then combine them in pairs and rearrange:
$$\langle \bar{B}_3-B_1,\bar{N}_2\rangle +\langle
B_3-\bar{B}_1,N_2\rangle =\langle B_1-\bar{B}_1,N_2\rangle
+\langle B_3-\bar{B}_3,N_2\rangle$$ (and cyclic). Adding up the
results, we find:
$$L(\theta)=\langle \bar{B}_1-B_1,N_1\rangle+\langle
\bar{B}_2-B_2,N_2\rangle +\langle \bar{B}_3-B_3,N_3\rangle,$$ so
we see that the total length of the network is independent of the
$x_i$, and in fact:
$$L(\theta)=\frac{\sqrt{3}}2L_{hex}(\theta).$$

\emph{5.2. Existence in a parallel curve.}

\begin{figure} \psset{unit=0.75cm}
\begin{center}
\SpecialCoor
\begin{pspicture}(-5,-4)(5,7)
\pnode(0,0){B1}\psline{->}([angle=300,nodesep=4.5]B1)([angle=120,nodesep=3]B1)
\pnode(6.1,7.8){B3bar}\psline{->}([angle=0,nodesep=0.5]B3bar)([angle=180,nodesep=4]B3bar)
\pnode(5.2,9){B2}\psline{->}([angle=60,nodesep=1]B2)([angle=240,nodesep=6]B2)
\pnode(-3.4,5.8){B1bar}\psline{->}([angle=120,nodesep=3]B1bar)([angle=300,nodesep=2.5]B1bar)
\pnode(-5.5,1.4){B3}\psline{->}([angle=180,nodesep=2]B3)([angle=0,nodesep=3.5]B3)
\pnode(-4.5,-1.3){B2bar}\psline{->}([angle=240,nodesep=2.5]B2bar)([angle=60,nodesep=3.5]B2bar)
\pnode([angle=300,nodesep=3]B1){B1new}
\pnode([angle=120,nodesep=1.5]B1bar){B1barnew}
\pnode([angle=180,nodesep=1.7]B3){B3new}
\pnode([angle=240,nodesep=2]B2bar){B2barnew}
\pnode([angle=180,nodesep=0.5]B3bar){B3barnew}
\psccurve[showpoints=true,linecolor=green](B1new)(B3barnew)(B2)(B1barnew)(B3new)(B2barnew)
\pnode([angle=120,nodesep=1.5]B1){H1}
\psline[linecolor=blue]{->}(H1)([angle=60,nodesep=7.5]H1)
\pnode([angle=60,nodesep=7.5]H1){H3bar}\psline[linecolor=blue]{->}(H3bar)([angle=300,nodesep=1.5]H3bar)
\pnode([angle=300,nodesep=1.5]H3bar){H2}\psline[linecolor=blue]{->}(H2)([angle=180,nodesep=7.5]H2)
\pnode([angle=180,nodesep=7.5]H2){H1bar}\psline[linecolor=blue]{->}(H1bar)([angle=240,nodesep=5.9]H1bar)
\pnode([angle=240,nodesep=5.9]H1bar){H3}\psline[linecolor=blue]{->}(H3)([angle=300,nodesep=3.8]H3)
\pnode([angle=300,nodesep=3.8]H3){H2bar}\psline[linecolor=blue]{->}(H2bar)([angle=0,nodesep=5.9]H2bar)
\uput[0](B1new){$B_1$}\uput[300]([angle=300,nodesep=4.5]B1){$n_1$}
\uput[45](B3barnew){$\bar{B_3}$}\uput[0]([angle=0,nodesep=0.5]B3bar){$\bar{n_3}$}
\uput[0](B2){$B_2$}\uput[0]([angle=240,nodesep=6]B2){$n_2$}
\uput[80](B1barnew){$\bar{B_1}$}\uput[120]([angle=120,nodesep=3]B1bar){$\bar{n_1}$}
\uput[240](B3new){$B_3$}\uput[180]([angle=180,nodesep=2]B3){$n_3$}
\uput[290](B2barnew){$\bar{B_2}$}\uput[240]([angle=240,nodesep=2.5]B2bar){$\bar{n_2}$}
\uput[270](H1){$H_1$}\uput[90](H3bar){$\bar{H_3}$}\uput[0](H2){$H_2$}\uput[70](H1bar){$\bar{H_1}$}
\uput{0.3}[35](H3){$H_3$}\uput[270](H2bar){$\bar{H_2}$}\uput[45]([angle=0,nodesep=5.9]H2bar){$H_1^*$}
\uput{0.5}[240]([angle=70,nodesep=3.3]B2bar){$R_3$}
\uput{0.2}[270]([angle=180,nodesep=1.4]B3bar){$\bar{R_3}$}
\uput{0.2}[90]([angle=70,nodesep=3.7]H1){(+)}
\uput[90]([angle=180,nodesep=3.7]H2){(+)}
\uput[70]([angle=300,nodesep=1]H3bar){(-)}
\end{pspicture}
 
\caption{Successive projection can lead to a `twisted chain'}
\end{center}
\end{figure}
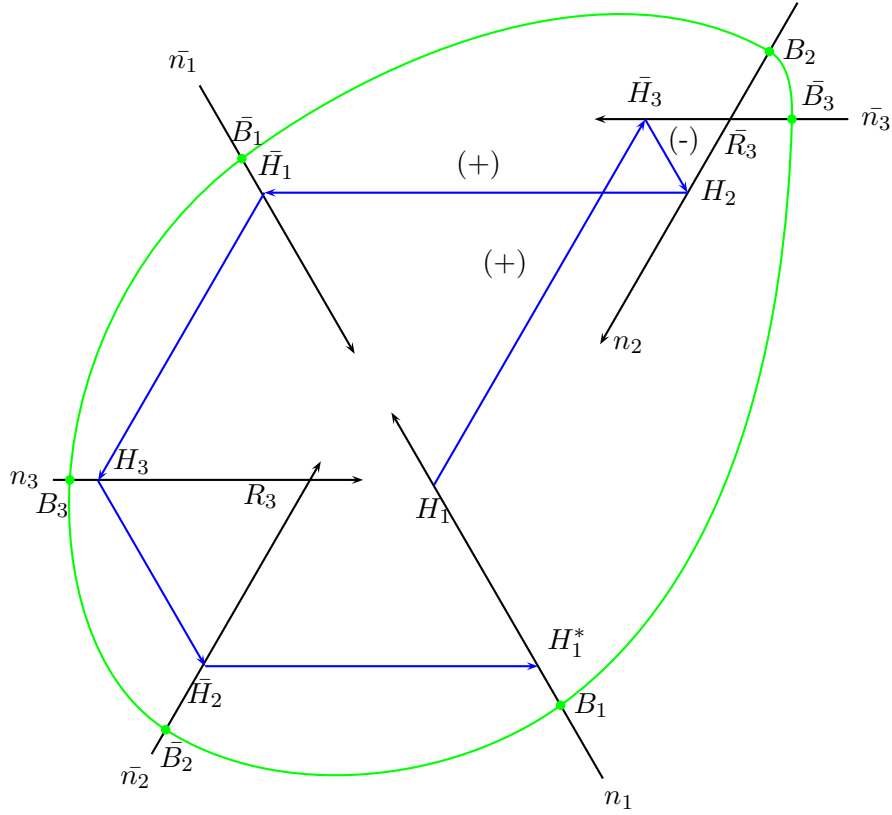
Vanishing of the holonomy is not the whole story. First, some
vertices of the cell could end up being outside the domain;
second, if we are not careful the `successive projection'
construction could lead to a `twisted cell'. (Fig.17; this example
has non-zero holonomy, but one sees easily that the problem may
also occur when the holonomy vanishes.)

The parameter $x=x_1$ of a hexagonal cell plays a secondary role
in existence considerations. Since consecutive interior vertices
of a hexagonal cell lie on consecutive normals, it is always
possible to `slide' the vertices along their normals until one
obtains a degenerate 5-vertex cell with one vertex at a `hexagonal
intersection' $R_I$ (where the interior angle is $\pi/3$) and four
vertices with interior angles $2\pi/3$: two `triangular
intersections' $P_I,P_J$ and two `parallelogram points' $U_I,V_J$.
After relabelling, we may assume the 5-cell is the chain (shown in
Fig. 15 in the case of nonzero holonomy):
$$V_3 \stackrel{pr_{3\bar{2}}}{\rightarrow} \bar{P_2}\stackrel{pr_{\bar{2}1}}{\rightarrow}R_1
\stackrel{pr_{1\bar{3}}}{\rightarrow}R_1\stackrel{pr_{\bar{3}2}}{\rightarrow}
P_1 \stackrel{pr_{2\bar{1}}}{\rightarrow}
\bar{U_1}\stackrel{pr_{\bar{1}3}}{\rightarrow} {V_3},$$ obtained
by successive projection along consecutive normals (taking $R_1$
to map to itself under $pr_{n_1\rightarrow\bar{n_3}}$). We also
used the fact that
$V_3:=pr_{\bar{2}3}(\bar{P_2})=pr_{\bar{1}3}(\bar{U}_1)$, which
follows from $hol(\theta)=0$. Existence of a 5-cell of this form
(entirely contained in $D$) is both necessary and sufficient for
the existence of a hexagonal cell in $D$. \vspace{.2cm}

Now, examination of Fig.17 shows the reason we get a `twisted'
cell is that $\bar{H}_3$ and $H_2$ are `ahead' of $\bar{R}_3$ (on
$\bar{n}_3$, resp. $n_2$); unlike, say, $H_3$ and $\bar{H}_2$
which are `behind' $R_3$ (on $n_3$, resp. $\bar{n}_2$). The normal
lines are oriented (by the inner unit normal vectors to $\cal C$),
so `ahead' and `behind' have meaning. Thus we need to make sure
the successive projections are always `behind' the $R_I$ along the
corresponding normals. \vspace{.2cm}

More precisely, in a `chain' obtained by projection on consecutive
normals, the angles between consecutive edges of the chain
($\pi/3$ or $2\pi/3$) depend on the relative \emph{orientations}
of the corresponding consecutive projections. Define the mapping
$pr_{n_I\rightarrow n_{I+1}}:X_I\mapsto X_{I+1}$ to be
\emph{positive} if $X_{I+1}-X_I$ has the  direction of $N_{I-1}$,
\emph{negative} if it has the direction of
$-N_{I-1}$.\vspace{.2cm}

This depends on considering three consecutive normals, which in
Fig.18 we denote by $n_1,\bar{n_3},n_2$. Three examples of
3-vertex chains are shown (from a point on $n_1$ to a point on
$n_2$). For $X_1\rightarrow \bar{X}_3\rightarrow X_2$, a (+)
projection is followed by a (-) projection, resulting in an angle
$\pi/3$ at $\bar{X}_3$; the reason $\bar{X}_3\rightarrow X_2$ is
(-) is that $\bar{X}_3$ and $X_2$ are both `ahead' of $\bar{R}_3$
on their respective normals. For
$X_1'\rightarrow\bar{X}_3'\rightarrow X_2'$, both $X_1'$ and
$\bar{X}_3'$ are `ahead' of $\bar{R}_3$; hence both projections
are (-), and the angle at $\bar{X}_3'$ is $2\pi/3$. Finally, with
$X_1'',\bar{X}_3''$ both `behind' $R_1$ and $\bar{X}_3''$, $X_2''$
`behind' $\bar{R}_3$, the chain $X_1''\rightarrow
\bar{X}_3''\rightarrow X_2''$ results from two (+) projections,
and the angle at $\bar{X}_3''$ is $2\pi/3$.\vspace{.2cm}

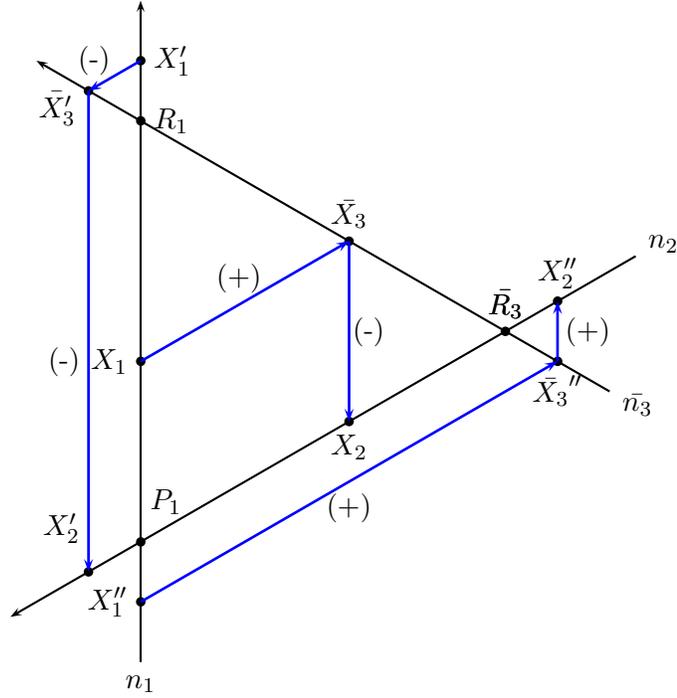
\begin{figure} \psset{unit=0.8cm}
\begin{center}
\SpecialCoor
\begin{pspicture}(-2,-1)(6,7)
\pnode(0,0){X1pp} \pnode([angle=30,nodesep=8]X1pp){X3pp}
\pnode([angle=90,nodesep=8]X1pp){R1}
\pnode([angle=90,nodesep=1]X1pp){P1}
\pnode([angle=90,nodesep=9]X1pp){X1p}
\pnode([angle=90,nodesep=4]X1pp){X1}
\pnode([angle=90,nodesep=10]X1pp){n1b}
\pnode([angle=270,nodesep=1]X1pp){n1a}
\pnode([angle=150,nodesep=1]X3pp){R3bar}
\pnode([angle=150,nodesep=4]X3pp){X3bar}
\pnode([angle=330,nodesep=1]X3pp){n3a}
\pnode([angle=150,nodesep=10]X3pp){n3b}
\pnode([angle=210,nodesep=1]P1){X2p}
\pnode([angle=30,nodesep=4]P1){X2}
\pnode([angle=30,nodesep=9.5]P1){n2a}
\pnode([angle=210,nodesep=2.5]P1){n2b}
\pnode([angle=30,nodesep=8]P1){X2pp}
\pnode([angle=90,nodesep=8]X2p){X3pbar}
\psline{->}(n3a)(n3b) \psline{->}(n1a)(n1b)\psline{->}(n2a)(n2b)
\psdots(X1pp)(X3pp)(R1)(P1)(X1p)(X1)(R3bar)(X3bar)(X2p)(X2)(X2pp)(X3pbar)
\psline[linecolor=blue,linewidth=1pt]{->}(X1pp)(X3pp)
\psline[linecolor=blue,linewidth=1pt]{->}(X3pp)(X2pp)
\psline[linecolor=blue,linewidth=1pt]{->}(X1)(X3bar)
\psline[linecolor=blue,linewidth=1pt]{->}(X3bar)(X2)
\psline[linecolor=blue,linewidth=1pt]{->}(X1p)(X3pbar)
\psline[linecolor=blue,linewidth=1pt]{->}(X3pbar)(X2p)
\uput[180](X1pp){$X_1''$}\uput{0.5}[60](P1){$P_1$}\uput[180](X1){$X_1$}\uput[0](R1){$R_1$}\uput[0](X1p){$X_1'$}
\uput[270](n1a){$n_1$}
\uput[30](n2a){$n_2$}\uput[90](X2pp){$X_2''$}\uput[90](R3bar){$\bar{R_3}$}\uput[270](X2){$X_2$}\uput{0.5}[120](X2p){$X_2'$}
\uput[90](R3bar){$\bar{R_3}$}\uput[330](n3a){$\bar{n_3}$}\uput[210](X3pbar){$\bar{X_3'}$}\uput[90](X3bar){$\bar{X_3}$}
\uput{0.5}[180](X1p){(-)} \uput[270](X3pp){$\bar{X_3}''$}
\uput{1}[180](X1){(-)}\uput{1.2}[270](X2){(+)}\uput{1}[0](R3bar){(+)}\uput{2}[180](R3bar){(-)}\uput{1.5}[200](X3bar){(+)}

\end{pspicture}
 
\caption{Orientations of consecutive projections on normals}
\end{center}
\end{figure}

A \emph{sufficient} condition for a 5-cell obtained by consecutive
projection from $R_1$ (backwards to $V_3$ and forward to
$\bar{U}_1$) to have the correct angles is that its vertices
always lie `behind' the intersections of normals in question; then
all projections will be (+), and all angles $2\pi/3$ (except at
$R_1$), as desired. (The condition is not necessary; having all
projections be (-) would achieve the same result.) This
requirement can be expressed in terms of the functions describing
the positions of points on the chain along their normals, as
follows.

$$\begin{array}{rcl}
R_1\leq \bar{R_3}(\bar{n}_3),\quad P_1\leq
\bar{R_3}(n_2)&\Leftrightarrow&
\bar{t}_3\leq \bar{s}_3\mbox{ and }v_2\leq t_2\\
P_1\leq R_2(n_2),\quad\bar{U}_1\leq
R_2(\bar{n}_1)&\Leftrightarrow& v_2\leq s_2\mbox{ and
}\bar{\sigma_1}\leq \bar{t}_1\\
\bar{P}_2\leq \bar{R}_2(\bar{n}_2),\quad R_1\leq
\bar{R}_2(n_1)&\Leftrightarrow& \bar{u}_2\leq \bar{s}_2\mbox{ and
}s_1\leq t_1\\
V_3\leq R_3(n_3),\quad \bar{P}_2\leq
R_3(\bar{n}_2)&\Leftrightarrow& \tau_3\leq s_3\mbox{ and
}\bar{u}_2\leq \bar{t}_2\\
\bar{U}_1\leq \bar{R}_1(\bar{n}_1),\quad V_3\leq
\bar{R}_1(n_3)&\Leftrightarrow& \bar{\sigma}_1\leq \bar{s}_1\mbox{
and }\tau_3\leq t_3
\end{array}$$
(equality is always allowed, since it just means the minimal
5-cell has `collapsed' to fewer than five vertices- the circle is
an extreme example). The two inequalities on the right in each of
the first four lines are equivalent to each other, even if
$hol(\theta)\neq 0$ (then we have an open 5-chain, from $V_3$ to
$\bar{U}_1$); if $hol(\theta)=0$, the same is true for the last
line.\vspace{.2cm}

We now use the previously computed expressions (5.6) for the
functions $u_I, v_I$, etc. appearing on the right-hand side to
express the inequalities in terms only of $\omega$ and $w$; this
yields six inequalities (two for the last line, since we don't
assume $hol(\theta)=0$ at this point.) To simplify them, we use
the easily verified identity:
$$w(\theta)+w(\theta+\frac{\pi}3)+w(\theta+\frac{2\pi}3)=
-\omega(\theta)-\omega(\theta+\frac{\pi}3).$$ The result is (in
the same order as above);
$$\begin{array}{lll}
&(A)&w(\theta+\frac{\pi}3)+\omega(\theta)\geq 0\\
&(B)&w(\theta)\leq 0\\
&(C)&w(\theta)+\omega(\theta+\frac{\pi}3)\leq 0\\
&(D)&w(\theta+\frac{\pi}3)\geq 0\\
&(E)& \omega(\theta)\geq 0\mbox{ and
}(F)\omega(\theta+\frac{\pi}3)\leq 0.
\end{array}$$
Note that (B),(D),(E),(F) is a `minimal set' (i.e., imply (A) and
(C)).

If we now use the condition
$$hol(\theta)=w(\theta)+w(\theta+\frac{\pi}3)+w(\theta+\frac{2\pi}3)=
-\omega(\theta)-\omega(\theta+\frac{\pi}3)=0,$$ (E) and (F) are
equivalent; so in this case (B),(D) and (E) are sufficient to
guarantee the projection construction based  at $R_1=R(\theta)$
yields a convex 5-cell with the correct angles.

In addition, rather than starting at $R_1$ (corresponding to
$\theta$), we could have started at any other $R_I$, corresponding
to the $\pi/3$-translation orbit of $\theta$. That is, we only
need (B),(D) and (E) to hold for some $\pi/3$ translate of
$\theta$. And it turns out (somewhat surprisingly, given the
experience with double triodes) that this is \emph{always} true
(assuming $hol(\theta)=0$), as verified in the combinatorial
proposition 5.4 below. We summarize the conclusion, bearing in
mind that only the configuration of normals plays a
role.\vspace{.2cm}

\textbf{Proposition 5.2.} \emph{Let
$n_1,\bar{n}_3,n_2,\bar{n}_1,n_3,\bar{n}_2$ be six oriented lines
in $\mathbb R^2$, with unit direction vectors between consecutive
lines differing by $\pi/3$ rotations; assume the configuration has
zero holonomy. Then it supports a one-parameter family of convex
equiangular hexagonal cells, with one vertex on each line.}
\vspace{.2cm}

\textbf{Corollary 5.3.} \emph{ For any strictly convex curve $\cal
C$ in $\mathbb R^2$ (possibly with corners), all sufficiently far
outer parallel curves (${\cal C}_d$ for $d\geq d_0\geq 0$) support
hexagonal cells.}\vspace{.2cm}

\textbf{Proposition 5.4.} \emph{Let $w$ and $\omega$ be $\pi$-periodic
and $2\pi/3$-periodic functions on $S$ (resp.) Suppose
$\theta_0\in S$ is a solution of the equations:
$$\omega(\theta_0)+\omega(\theta_0+\frac{\pi}3)=0,\quad
w(\theta_0)+w(\theta_0+\frac{\pi}3)+w(\theta_0+\frac{2\pi}3)=0.$$
Then some $\pi/3$ translate $\hat{\theta}$ of $\theta_0$
satisfies, in addition to these two equations, also the
inequalities:
$$w(\hat{\theta})\leq 0,\quad w(\hat{\theta}+\frac{\pi}3)\geq
0,\quad \omega(\hat{\theta})\geq 0.\quad (5.7)$$}

\emph{Proof.} (By exhaustive listing of cases.) The sign of
$\omega(\theta_0)$ determines that of $\omega(\theta)$ for all
$\pi/3$-translates $\theta$ of $\theta_0$, given $2\pi/3$
periodicity and the first equation. For $w$, given the sign of
$w(\theta_0)$ the second equation implies 3 possible sign
combinations for $w(\theta_0+\pi/3)$ and $w(\theta_0+2\pi/3)$,
which then determine the signs at the other translates. This gives
12 cases, and in each case we can find translates satisfying all
three inequalities desired.

We proceed to list the 12 cases. On each half-line, the four signs
are those of $\omega(\theta_0),w(\theta_0), w(\theta_0+\pi/3),
w(\theta_0+2\pi/3)$ (in this order); the last entry is the
translate $\hat{\theta}$ of $\theta_0$ satisfying the sign
conditions in the proposition.

$$\begin{array}{cccclccccl}
+&-&+&+&\theta_0&+&+&-&-&\theta_0+2\pi/3\\
+&-&-&+&\theta_0+4\pi/3&+&+&-&+&\theta_0+4\pi/3\\
+&-&+&-&\theta_0&+&+&+&-&\theta_0+2\pi/3\\
-&-&+&+&\theta_0+\pi&-&+&-&-&\theta_0+5\pi/3\\
-&-&-&+&\theta_0+\pi/3&-&+&-&+&\theta_0+\pi/3\\
-&-&+&-&\theta_0+\pi&-&+&+&-&\theta_0+5\pi/3
\end{array}$$

\vspace{.4cm}

\emph{5.3 Existence under curvature conditions.} In this
subsection we show that sufficiently strong `pinching conditions'
on the radius of curvature imply all 24 points of the
configuration are inside the domain; combined with the conclusion
of the previous sections, this shows hexagonal cells exist in this
case.

The points of the configuration on the normal $n_1=n(\theta)$ are:
$$P_1, P_2, R_1, \bar{R}_2,U_1,V_1,$$
and all can be written in the form $B_1+(2/\sqrt{3})f_IN_1$, where
$f_I$ is, respectively:
$$u(\theta),v(\theta),s(\theta),t(\theta),\sigma(\theta),\tau(\theta).$$
Recall from section 3 the chord $n_1\cap D$ has endpoints
$B(\theta), B(\theta_*)$ (with $\theta_*\in (\theta,
\theta+2\pi)$) and length:
$$d(\theta)=\langle B(\theta_*)-B(\theta),N(\theta)\rangle.$$
Thus, we seek conditions that imply
$0<f_I(\theta)<(\sqrt{3}/2)d(\theta)$, for each $f_I$ given above
and all $\theta$.

For $u$ and $v$, this was done in Proposition 3.2: we showed that
$r_{max}<(1+\sqrt{3}/2)r_{min}$ implies $\theta_*\in
[\theta+2\pi/3,\theta+4\pi/3]$ for all $\theta$, and that this
condition implies both $u(\theta)$ and $v(\theta)$ are in
$(0,(\sqrt{3}/2)d(\theta))$. \vspace{.2cm}

Our first observation is that $s(\theta)$ and $t(\theta)$ are
always \emph{positive}:
$$\begin{array}{rcl}
s(\theta)=\langle
B(\theta+\frac{\pi}3)-B(\theta),T(\theta+\frac{\pi}3)\rangle
&=&\int_0^{\pi/3}r(\tau+\theta)\cos(\tau-\pi/3)d\tau>0,\\
t(\theta)=\langle
B(\theta)-B(\theta-\frac{\pi}3),T(\theta-\frac{\pi}3)\rangle
&=&\int_{-\pi/3}^0r(\tau+\theta)\cos(\tau+\pi/3)d\tau>0.
\end{array}\quad(5.8)$$

\emph{Remark 5.5.} From relations (5.2),(5.3), we see that this
implies $v(\theta)>w(\theta+\pi/3)$ and
$u(\theta)>-w(\theta+2\pi/3)$ for each $\theta$, and we conclude:
(i) curves of constant width ($w\equiv 0$) always support triodes,
and (ii) at a zero of $hol$ satisfying inequalities (5.7), we
automatically have $v(\theta)>0$ (but not necessarily
$u(\theta)>0$).\vspace{.2cm}

The following two lemmas deal with the pairs $(s,t)$ and
$(\sigma,\tau)$.\vspace{.2cm}

\textbf{Lemma 5.6} \emph{Assume $\theta^*\in
[\theta+\frac{2\pi}3,\theta+\frac{4\pi}3]$, for all $\theta$. Then
if $r_{max}<(4/3)r_{min}$, we have $s(\theta),t(\theta)$ both
(positive and) less than $(\sqrt{3}/2)d(\theta)$, for all
$\theta$.}

\emph{Proof.} (i) From (5.8) we see that, to show
$s<(\sqrt{3}/2)d$, we need:
$$\frac 12 \int_0^{\pi/3}r(\tau+\theta)\cos \tau
d\tau+\frac{\sqrt{3}}2\int_0^{\pi/3}r(\tau+\theta)\sin \tau
d\tau<\frac{\sqrt{3}}2\int_0^{\theta_*-\theta}r(\tau+\theta)\sin
\tau d\tau,$$ or:
$$\int_0^{\pi/3}r(\tau+\theta)\cos \tau
d\tau<\sqrt{3}[\int_{\pi/3}^{2\pi/3}r(\tau+\theta)\sin \tau d\tau
+\int_{2\pi/3}^{\theta_*-\theta}r(\tau+\theta)\sin \tau d\tau].$$
We estimate the left-hand side from above by
$\frac{\sqrt{3}}2r_{max}$. If $\theta_*-\theta\leq \pi$, estimate
the right-hand side from below by $\sqrt{3}r_{min}$ (using only
the first integral on the right). If $\theta_*-\theta\geq \pi$,
the right-hand side is bounded below by:
$$\sqrt{3}[r_{min}+r_{min}\int_{2\pi/3}^{\pi}\sin \tau d\tau
-r_{max}\int_{\pi}^{4\pi/3}|\sin \tau|d\tau]=\sqrt{3}[\frac 32
r_{min}-\frac 12r_{max}].$$ This gives the condition
$r_{max}<(3/2)r_{min}$.

(ii) To show $t<(\sqrt{3}/2)d$, from (5.8) we need the condition:
$$\begin{array}{rcl}\int_{-\pi/3}^0r(\tau+\theta)\cos \tau
d\tau&+&\sqrt{3}\int_{-\pi/3}^0r(\tau+\theta)|\sin
\tau|d\tau\\&<&\sqrt{3}[\int_0^{2\pi/3}r(\tau+\theta)\sin \tau
d\tau+ \int_{2\pi/3}^{\theta_*-\theta}r(\tau+\theta)\sin \tau
d\tau.]\end{array}$$ One verifies easily that the left-hand side
is bounded above by $\sqrt{3}r_{max}$. For the right-hand side,
there are two cases: if $\theta_*\in (\frac{2\pi}3,\pi]$, the
right-hand side is greater than
$\sqrt{3}r_{min}[1-\cos(\theta_*-\theta)]>(3\sqrt{3}/2)r_{min}.$
If $\theta_*-\theta\in [\pi,4\pi/3)$, the right-hand side is
bounded below by:
$$\sqrt{3}[r_{min}\int_0^{\pi}\sin \tau
d\tau-r_{max}\int_{\pi}^{4\pi/3}|\sin \tau
|d\tau]=\sqrt{3}[2r_{min}-\frac 12 r_{max}].$$ This gives the
condition $(3\sqrt{3}/2)r_{max}<2\sqrt{3}r_{min}$, or
$r_{max}<(4/3)r_{min}$, which is more stringent than that in part
(i) above.\vspace{.2cm}

\textbf{Lemma 5.7.}\emph{ Assume $\theta_*\in
[\theta+2\pi/3,\theta+4\pi/3]$, for all $\theta$. Then if
$(4+3\sqrt{3})r_{max}<(4+4\sqrt{3})r_{min}$, we have:
$0<\tau<(\sqrt{3}/2)d$ and $0<\sigma<(\sqrt{3}/2)d$.
}\vspace{.2cm}

\emph{Proof.} Using $\tau(\theta)=s(\theta)-w(\theta)$, to prove
the statement for $\tau$ we need to verify:
$$\begin{array}{rcl}
0<\frac 12\int_0^{\pi/3}r(\tau+\theta)\cos \tau
d\tau&+&\frac{\sqrt{3}}2\int_0^{\pi/3}r(\tau+\theta)\sin \tau
d\tau\\&-&\int_0^{\pi/2}r(\tau+\theta)\cos \tau
d\tau+\int_{\pi/2}^{\pi}r(\tau+\theta)|\cos \tau|d\tau\\
&<&\frac{\sqrt{3}}2(\int_0^{2\pi/3}r(\tau+\theta)\sin \tau d\tau
+\int_{2\pi/3}^{\theta_*-\theta}r(\tau+\theta)\sin \tau
d\tau).\end{array}$$

For $\tau>0$, we see it is enough to show that:
$$0<\frac{\sqrt{3}}4r_{min}+\frac{\sqrt{3}}4r_{min}-r_{max}+r_{min},$$
so the condition is $r_{max}<(1+\sqrt{3}/2)r_{min}$. For the upper
bound on $\tau$, we again consider two cases. If $\theta_*\in
[\theta+2\pi/3,\theta+\pi]$, we need:
$$\frac{\sqrt{3}}4r_{max}+\frac{\sqrt{3}}4r_{max}-r_{min}+r_{max}<\frac{\sqrt{3}}2
[\frac 32 r_{min}+r_{min}|\cos(\theta_*-\theta)+\frac 12|],$$ or
$(4+2\sqrt{3})r_{max}<(4+3\sqrt{3})r_{min}$.

If $\theta_*\in [\theta+\pi,\theta+4\pi/3]$, we need:
$$\frac{\sqrt{3}}2r_{max}-r_{min}+r_{max}<\frac{\sqrt{3}}2[\frac
32 r_{min}+r_{min}\int_{2\pi/3}^{\pi}\sin
\tau-r_{max}\int_{\pi}^{4\pi/3}|\sin\tau|];$$ this gives the
condition $(4+3\sqrt{3})r_{max}<(4+4\sqrt{3})r_{min}$, the more
stringent of the three.

The proof for $\sigma$ is almost identical (and gives the same
constant), so we omit it. \vspace{.2cm}

We summarize the discussion in the following proposition.

\textbf{Proposition 5.8.} \emph{Let $\cal C$ be a strictly convex
curve ($p\in {\cal P}_{pw}$) satisfying the radius of curvature
bounds
$\frac{r_{max}}{r_{min}}<\frac{4+4\sqrt{3}}{4+3\sqrt{3}}\sim
1.188\ldots$ Then for each $\theta\in S$, all 24 points of the
`hexagonal configuration' are inside $\cal C$. Thus $\cal C$
(stably) supports at least two geometrically distinct `bands' of
hexagonal cells (corresponding to the global max and min of the
`hexagonal perimeter' $L_{hex}$).}\vspace{.3cm}

\emph{5.4. Examples of nonexistence.}

If $\cal C$ supports a hexagonal cell, in particular it must
support a minimal 5-cell as described in section 5.2: for some
`hexagonal intersection' $R_1=R(\theta)$,the 5-cell has the form
$V_3\bar{P}_2R_1P_1\bar{U}_1$, where the chain is obtained by
backwards and forward projection on consecutive normals, beginning
at $R_1$; in particular, $P_1$ and $\bar{P}_2$ are `triangle
intersections'. All the points in the 5-chain must be inside $\cal
C$. Thus we have the \emph{necessary condition}: if $\theta$ with
$hol(\theta)=0$ corresponds to a hexagonal cell, then for some
$\hat{\theta}$ among the $\pi/3$-translates of $\theta$ we must
have $u(\hat{\theta})>0$ (so $P_1(\hat{\theta})\in D$) and
$u(\hat{\theta}-\pi/3)>0$ (so $\bar{P}_2(\hat{\theta})\in D$).
\vspace{.2cm}

Suppose $\omega=u-v$ is $\pi$-periodic, and therefore in fact
$\pi/3$-periodic. Then $hol\equiv 2\omega$, so at a zero of $hol$
we have $P_1=P_2=P_3:=P$ and
$\bar{P_1}=\bar{P_2}=\bar{P_3}:=\bar{P}$. If, in addition, we find
that $u(\theta)<0$ and $u(\theta+\pi)<0$ for \emph{every} zero of
$\omega$ (or $hol$) in some $\pi/3$ fundamental domain, then at
every such zero both $P$ and $\bar P$ are outside the domain; thus
$\cal C$ does not support a hexagonal cell (or a triode). This
leads to the following construction. \vspace{.2cm}

\textbf{Proposition 5.9.} \emph{Assume $p_0$ has the properties:
(i) strict convexity ($p_0\in {\cal P}_{smooth}$); (ii)constant
height ($\omega_0\equiv 0$); (iii)$u_0<0$ in some interval
$I=(-\delta,\delta)$, and also in
$I+\pi=(-\delta+\pi,\delta+\pi)$, for some $0<\delta<\pi/6$. Let
$\tilde{p}$ have the properties: (iv)$\tilde{p}$ is
$\pi/3$-periodic; (v) all the zeros of $\tilde{\omega}$ in a
$\pi/3$-fundamental domain (say, $[-\frac{\pi}6,\frac{\pi}6]$) are
contained in $I$. Then:
$$p:=\tilde{p}+\lambda p_0$$ has the properties, for $\lambda$
sufficiently large: (a) strict convexity; (b) at each zero of
$hol$ (equivalently, of $\omega$, since $\omega=\tilde{\omega}$ is
$\pi/3$-periodic), $P$ and $\bar{P}$ are outside $D$. In
particular, $\cal C$ supports no hexagonal cells, and no triodes.
(c) Any support function sufficiently $C^2$-close to $p$ also
satisfies (a) and (b).}

\emph{Proof.} If $\lambda$ is large enough, the zeros of
$\omega=\tilde{\omega}$ in a fundamental domain containing
$\theta=0$ are all in $I$ (where $u<0$ for $\lambda$ large), so
for any $\theta$ with zero holonomy in this fund. domain we have
$P\not \in D$; while the zeros of $\omega$ (or $hold$ in a fund.
domain containing $\pi$ are all in $I+\pi$, where again $u<0$, so
$\bar{P}\not \in D$.  From the discussion above, (b) follows. (a)
and (c) are clear. \vspace{.3cm}

\emph{Example 5.1} To get an explicit example, we modify Example
3.3; take for $p_0$ the support function of a constant height
$\epsilon$-biangle, for $\epsilon>0$ sufficiently small, or a
smoothed version, such as $p_{\epsilon}$ in example 3.3; then
$I=(-0.34,0.34)$. With $m=\sqrt{2}/2$, set:
$$q(\theta)=\cos(6\theta)/(1-m\sin(6\theta)).$$ Then
 $q$ is $\pi/3$-periodic, with critical points
$\pi/24$ and $\pi/8$ in $[0,\pi/3)$. Defining
$\tilde{p}=q(\theta+\pi/12)$, we find $\tilde{\omega}=3\tilde{p}'$
has zeros at $\pm \pi/24\sim \pm 0.1309$ in the fundamental domain
$[-\pi/12,\pi/4)$, both in $I$- so $\tilde p$ satisfies (iv) and
(v) of the proposition.\vspace{.2cm}

\textbf{Corollary 5.10.} \emph{If $p_0\in {\cal P}_{smooth}$ has
constant height and satisfies $u_0<0$ on some interval $I$
\emph{and} on its translate $I+\pi$, then arbitrarily $C^2$-close
to $p_0$ one finds strictly convex curves $p\in {\cal P}_{smooth}$
with the property that curves in an open $C^2$-neighborhood of $p$
support no hexagonal cells or triodes.}\vspace{.2cm}

\emph{Remark 5.6- multiple non-existence.} It is natural to try to
refine this construction, so as to obtain domains supporting no
triodes, double triodes, or hexagonal cells. Unfortunately our
examples of non-existence for double triodes (section 4) rely on
control of the derived width function $w(\theta)$, while the class
of convex curves of \emph{constant height} with $u<0$ somewhere
(on which the examples without triodes are based) appears to be
too small to allow fine control of $w$. While the construction of
section 4 may be used to find smooth convex curves of constant
height without double triodes, all known examples have $u>0$
everywhere (Remark 4.4). On the other hand, there seems to be no
fundamental reason why a convex domain without any Steiner
networks would be an impossibility.



\begin{thebibliography}{99}
\bibitem[Bronsard-Reitich]{Bronsard-Reitich} Bronsard,L., Reitich,
F., \emph{On three-phase boundary motion and the singular limit of
a vector-valued Ginzburg-Landau equation}  Arch. Rational Mech.
Anal. \textbf{124} (1993),  no. 4, 355--379.

\bibitem[Ivanov-Tuzhilin]{IvanovTuzhilin} Ivanov, A., Tuzhilin, A.
 \emph{Branching solutions to one-dimensional variational problems.}
 World Scientific Publishing Co., Inc., River Edge, NJ, 2001. xxii+342 pp.
 ISBN: 981-02-4060-0

\bibitem[Mantegazza et al.]{Mantegazza et al.} Mantegazza,C.,
Novaga,M., Tortorelli,V., \emph{Motion by curvature of planar
networks} Ann. Sc. Norm. Super. Pisa Cl. Sci. \textbf{5}(3)
(2004). no. 2, 235--324.

\bibitem[Mese-Yamada]{MeseYamada}Mese, C., Yamada, S.
 \emph{The parameterized Steiner problem and the singular plateau Problem via energy.}
 Trans. Amer. Math. Soc.  \textbf{358}  (2006), no.7, 2875--2895

\bibitem[Tabachnikov]{Tabachnikov} S.Tabachnikov, \emph{The
four-vertex theorem revisited}, American Mathematical Monthly
\textbf{102}(10) (1995), 912--916.

\bibitem[Yaglom-Boltianskii]{Yaglom-Boltianskii} I.M. Yaglom,
V.G.Boltianski, \emph{Convex Figures}, Holt, Rinehart and Winston,
1961.

\bibitem[Ikota-Yanagida]{Ikota-Yanagida}Ikota, R., Yanagida, E.\emph{Stability of stationary interfaces
of binary-tree type} Calc. Var. Partial Differential Equations
\textbf{22} (2005), no. 4, 375--389.
\end{thebibliography}
\end{document}